\documentclass[a4paper,fleqn]{cas-sc}

\usepackage{longtable}
\usepackage{enumitem}
\usepackage{appendix}
\usepackage{bm}
\usepackage{color}
\usepackage{amsmath}
\usepackage{amsthm}
\usepackage{float}
\usepackage{braket}
\usepackage{algorithmic}
\usepackage{algorithm}
\usepackage[authoryear]{natbib}

\def\tsc#1{\csdef{#1}{\textsc{\lowercase{#1}}\xspace}}
\tsc{WGM}
\tsc{QE}
\tsc{EP}
\tsc{PMS}
\tsc{BEC}
\tsc{DE}

\begin{document}
\let\WriteBookmarks\relax
\let\printorcid\relax 

\shorttitle{\rmfamily\small Collaborating Unmanned Aerial Vehicle and Ground Sensors for Urban Signalized Network Traffic Monitoring} 

\shortauthors{Jiarong Yao et al.}

\title[mode = title]{Collaborating Unmanned Aerial Vehicle and Ground Sensors for Urban Signalized Network Traffic Monitoring}

\author[1]{Jiarong Yao}[orcid=0000-0003-3058-189X] 
\ead{joanna.yao@polyu.edu.hk} 

\author[2]{Chaopeng Tan}[orcid=0000-0003-4737-5304] 
\ead{chaopeng.tan@tu-dresden.de} 
\cormark[1]

\author[2]{Meng Wang}

\author[1]{Wei Ma}

\address[1]{Department of Civil and Environmental Engineering, The Hong Kong Polytechnic University, Hung Hom Kowloon, Hong Kong}
\address[2]{Chair of Traffic Process Automation, Technische Universität Dresden, Dresden, Germany}

\cortext[1]{Corresponding author} 

\begin{abstract}
Reliable estimation of network-wide traffic states is essential for urban traffic management. Unmanned Aerial Vehicles (UAVs), with their airborne full-sample continuous trajectory observation, bring new opportunities for traffic state estimation. In this study, we will explore the optimal UAV deployment problem in road networks in conjunction with ground sensors, including connected vehicle (CV) and loop detectors, to achieve more reliable estimation of vehicle path reconstruction as well as movement-based arrival rates and queue lengths. Oriented towards reliable estimation of traffic states, we propose an index, feasible domain size, as the uncertainty measurement, and transform the optimal UAV deployment problem into minimizing the observation uncertainty of network-wide traffic states. Given the large-scale and nonlinear nature of the problem, an improved quantum genetic algorithm (IQGA) that integrates two customized operators is proposed to enhance neighbor searching and solution refinement, thereby improving the observability of UAV pairs. Evaluation was conducted on an empirical network with 18 intersections. Results demonstrated that a UAV fleet size of 7 is sufficient for traffic monitoring, with more than 60\% of network-wide observation uncertainty reduced. Through horizontal comparison with three baselines, the optimal UAV location scheme obtained by the proposed method can reach an improvement of up to 7.23\% and 5.02\% in the estimation accuracy of arrival rate and queue length, respectively. The proposed IQGA is also shown to be faster in solution convergence than the classic QGA by about 9.22\% with better exploration ability in optimum searching.
\end{abstract}


\maketitle

\section{Introduction}

Reliable and comprehensive traffic state monitoring plays a critical role in effective urban traffic management. Existing studies have demonstrated that more accurate traffic state estimates can reduce 12\%-20\% average vehicle delay at intersections \citep{tan2025robust}, which not only improves travel efficiency but also effectively reduces fuel consumption and carbon emissions, contributing to sustainable urban development \citep{karkouch2016data}. 
Conventional traffic monitoring for urban signalized network has predominantly focused on two scales: network level and intersection level. The former involves estimating origin-destination (OD) or path flows \citep{tang2025}, while the latter focuses on capturing intersection dynamics such as vehicle arrivals and queue length \citep{Wu31122024}. Network-scale monitoring is essential for understanding the overall travel demand distribution and network-wise performance, yet it poses significant challenges in data collection and estimation due to the limitations of stationary sensors and the inherent difficulties in capturing complete travel patterns \citep{krishnakumari2020data, sun2024stochastic}. At the intersection level, obtaining accurate data on intersection traffic flows and queue lengths is equally critical, as such data informs adaptive signal control strategies and congestion management; however, traditional stationary sensors are often hampered by high installation costs, limited spatial coverage, and insufficient temporal resolution \citep{luo2025probabilistic, tan2021cumulative}. Although advancements in connected vehicle (CV) technology have provided a promising alternative by harnessing real-time vehicular data for traffic state estimation, e.g., OD matrix \citep{cao2021day}, path flow \citep{chen2021dynamic} estimation at the network level, and queue length \citep{tan2019cycle, hao2014cycle}, volume \citep{zheng2017estimating, yao2019sampled} estimation at intersections, they are fundamentally sampled observations and CV-based methods are thus subject to uncertainties and biases, particularly under conditions of low penetration rates. Therefore, these traditional ground sensors inherently lack the capability to provide comprehensive spatiotemporal observations of traffic flows.

Simultaneously, the rapid development of unmanned aerial vehicles (UAVs) has introduced an innovative perspective to the domain of traffic monitoring \citep{PONSPRATS2022102868, Prathibaa2023, Telikani2024, ge2025multi, lyu2025urban}. UAVs have the unique capability of continuously capturing high-resolution imagery and video over specific areas, thereby enabling the tracking of vehicle trajectories and real-time assessment of traffic conditions from a top-down view \citep{bisio2022systematic}. This aerial approach offers a distinct advantage by covering the full sample and providing a more complete picture of vehicular movements. Ideally, if we have a fleet of UAVs hovering above all intersections during a time-of-day (TOD) period for monitoring, both the network-level and intersection-level traffic states can be perfectly obtained \citep{barmpounakis2020new}. However, such a practice is far from feasible as full coverage for all the signalized intersections may be costly regarding device procurement, configuration, and maintenance. Moreover, a 100\% coverage is sometimes redundant or unnecessary, considering the traffic flow correlation among intersections.

Given these considerations, this study uses limited UAV resources as a complement to network-wide traffic monitoring, aiming to collaborate ground sensors with UAVs for more robust and reliable traffic state estimation. Essentially, the research problem is to assign the locations of the limited number of UAVs to maximize the "benefits" of network traffic monitoring. Existing studies on sensor location optimization have primarily measured such benefits using path or flow coverage rates \citep{GENTILI2012227}, while some studies have employed entropy-based measures of uncertainty to assess the benefits \citep{tang2025novel}. However, these metrics fail to provide a unified framework for evaluating the benefits of heterogeneous sensors, such as UAVs and ground sensors (loop detectors and CVs) in this study, across both network- and intersection-level traffic monitoring dimensions. Therefore, this study proposes a new metric, the feasible domain size, to uniformly measure uncertainty of different traffic states under heterogeneous sensors. This uncertainty measure is then used to construct a UAV location optimization framework, providing exact intersection sets to maximize the benefits of the road network monitoring. 

The major contribution of this study is threefold: 1) As far as we know, this is the first study on collaborating ground sensors and UAVs for network-wide traffic state monitoring. 2) We propose a new metric, feasible domain size, to measure the uncertainty of traffic states, which provides a uniform framework to quantify the benefits of heterogeneous traffic sensors as well as traffic state indicators of different granularities. 3) We propose an uncertainty minimization framework that transforms the traffic state uncertainty optimization problem into a UAV location optimization problem, which improves the data utilization efficiency at the planning level and facilitates various traffic state estimation methods.

\section{Literature review}
\subsection{Detector location optimization}
Early work treats detector placement as a network sensor location problem (NSLP) for link-flow and origin-destination (OD) inference, relying mainly on detectors without ID information, e.g., loop detectors. \cite{yang1998optimal} proposed four canonical rules (OD coverage, maximum flow fraction/intercept, link independence) and integer programs; this set the agenda for count-based placement for OD estimation. \cite{ehlert2006optimisation} formalized the network count location problem (NCLP), selecting “informative” links under budget. \cite{bianco2006combinatorial} contributed a combinatorial analysis, proving NP-completeness in general and characterizing graph classes where polynomial solutions exist. Together, these papers established observability-driven and coverage-style formulations solved via ILP/heuristics. 

Later, the focus shifted to complete or partial link-flow observability and computational scalability. \cite{hu2009identification} gave constructive conditions for link-based applications; \cite{ng2012synergistic} introduced a node-based “synergistic” approach that avoids path enumeration; \cite{he2013graphical} provided a graphical method; and robust models now explicitly propagate measurement error \citep{xu2016robust}, with redundancy/failure-tolerant designs later explored and extended, e.g., \cite{salari2019optimization}; \cite{zhu2022network}. Recent theory shows submodularity can yield greedy algorithms with guarantees \citep{li2023submodularity}, while dynamic-state observability has been studied for density estimation in ODE/first-order traffic models \cite{hu2024sensor}. Key objectives include minimizing sensors subject to full observability, minimizing estimation variance/uncertainty, and ensuring robustness to failure under integer programming or greedy/submodular frameworks.

Recent advancements in detection techniques have facilitated the deployment of detectors like Automated Vehicle Identification (AVI) sensors, which provide vehicle re-identification between pairs (or sequences). A major stream of AVI deployment works on the recovery of OD demand, path flows, and link travel times, which treats sensor location as an observability problem under practical constraints. Early AVI-specific formulations used prior OD information to steer deployment and estimate both the OD matrix and routes; for instance, \cite{minguez2010optimal} cast “plate scanning” placement as an optimization guided by a priori OD distributions and showed how AVI reads enable joint OD and route estimation. \cite{castillo2008trip, castillo2010optimal} studied partial path observability with AVI (and link counts), including the effects of missed detections and how to use the order and timing of plate reads to improve inference. Building on these, \cite{cerrone2015vehicle} explicitly modeled the sequence of AVI placements to strengthen path differentiation; \cite{hadavi2019vehicle} then proposed fast greedy/meta-heuristic algorithms that scale to large networks while preserving route-flow observability guarantees. The generation of the candidate path set itself has a major impact: \cite{rinaldi2017exact} quantified how path‐set construction affects the solvability of the flow-observability problem and classified paths by information redundancy. More recently, \cite{sun2022reliable} addressed reliability by locating AVI sensors so OD demands remain uniquely identifiable even under sensor failures, avoiding an intermediate route-flow step and minimizing sensors given a target reliability. For joint OD and travel time estimation, multi-type/heterogeneous deployments have been studied. \cite{fu2022optimization} formulated a multi-type sensor location problem that simultaneously estimates OD and link travel times, using Kullback–Leibler divergence and explicitly accounting for covariance between measurements and detection errors; their model also incorporates deployment costs. \cite{zhu2023optimal} tackled link-level travel time coverage and accuracy with an AVI placement method (framed as a potential/benefit-increment game in follow-ups) and gave practical guidance on how many readers are needed and where. Separately, \cite{alvarez2022iterative} evaluated AVI deployments under uncertainty and popularized RMARE as a quality metric for flow estimation—often cited in subsequent deployment studies.

A second thread targets reconstructing vehicle paths (and derivative quantities) from AVI, often exploiting micro-level IDs and timestamps. \cite{sanchez2011dealing, sanchez2020new} showed that with plate reads (possibly sparse and noisy), one can partially infer paths and update OD estimates via Bayesian networks, and they developed placement rules that maximize the inferential value of ordered AVI observations while addressing error recovery. This line motivates sensor layouts that maximize the discriminability of competing paths. On the deployment side, \cite{fu2017stochastic} introduced a path-reconstruction-oriented sensor location model under uncertainty; later work considers mixing AVI with loop detectors to reduce cost while keeping reconstruction fidelity. Pushing the objective explicitly toward reconstruction, \cite{tang2025novel} proposed Path Reconstruction Entropy (PRE) and an AVI location model that minimizes PRE, prioritizing placements that most reduce ambiguity in individual paths.

In summary, it is evident that existing studies on detector deployment primarily focus on the benefits of network-level traffic state observations at the road network level, such as OD flows and path reconstruction, while neglecting the benefits of intersection-level traffic state observations. These intersection-level traffic states are crucial for evaluating and optimizing traffic signal systems. 

\subsection{UAV-based traffic monitoring}

Thanks to the advancements in vision algorithms and image processing, applications of UAVs in the transportation field especially in traffic monitoring is emerging \citep{outay2020applications}. Object detection and tracking is found as a key element to extract key information of road users from the field-of-view (FoV) from snapshots or videos collected by UAVs before further traffic flow analysis and traffic state assessment \citep{gohari2022involvement}. Thus, early studies using UAV for traffic monitoring mainly focus on the refinement of computer vision-related algorithms to guarantee the stability and accuracy of objects or trajectories extracted from UAV data \citep{bisio2022systematic}. The output traffic state indicators of such studies are simply derivatives of single trajectory like position, speed \citep{garcia2023trajectory} or deduced based on the classic traffic flow theory, such as flow, density, shockwave velocity \citep{wang2016improved, ke2016real}, etc. 

Different from other existing data sources for traffic monitoring, UAV detector has its uniqueness and irreplaceability. For UAV hovering detection mode, a UAV detector is able to capture the full sample of traffic flow within a certain temporal-spatial scale, which provides both macroscopic traffic flow indicators of flow, density and speed, as well as microscopic driving behaviors like lane-changing and car-following of every trajectory in FoV \citep{butilua2022urban}. As for UAV cruising detection mode, it is possible for a UAV detector to track the temporal-spatial evolution of traffic flows in the network scale, especially under a multi-UAV collaboration pattern \citep{gupta2021advances}. Compared with mobile sensors like connected vehicles which only provide sampled trajectories and fixed detectors which are inflexible to adapt to changing traffic demand patterns, UAV data realize the consistency between macro- and micro- properties, mobile sensing and static detection, which possess potential to create new paradigms for traffic monitoring and management.

Taking this one step further, researchers start to make more use of the UAV data for informative traffic monitoring, which can be shown in existing studies regarding dimensions in application scenario and analysis granularity. As a continuum flow facility, expressway is mostly seen in early literature as the research scenario, with the traffic monitoring objective of speed or density collection, trajectory reconstruction, or later macroscopic fundamental diagram (MFD) or shockwave for congestion monitoring \citep{li2021domain, bayraktar2023traffic}. Starting from fixed-point collection, \cite{yahia2022unmanned} utilized the congestion level of different freeway sections to back-nurture the UAV location for traffic monitoring, which is developed into a joint optimization of single UAV routing and traffic state estimation. Further, signalized roadway network scenario, which is more complicated than expressway scenario, has been studied, starting from isolated link-based or node-based UAV detection \citep{kaufmann2018aerial,barmpounakis2019accurate}. With the interference of signal controllers, traffic flows passing the intersection are interrupted due to difference in demand distribution and signal timing, thus platoon matching and platoon clustering are adopted by some researchers for path-based or movement-based traffic monitoring providing arrival volume or OD flow estimation \citep{coifman2006roadway, ke2016real, ke2018real}. 

With multiple UAVs available, more traffic state indicators can be obtained through the joint-monitoring of the UAV fleet within a network, like travel time, abnormal vehicle behavior, and congestion propagation, etc. \citep{zhu2018urban, barmpounakis2020new}, especially when they are linked up with downstream applications like signal control \citep{yao2025incorporating}, vehicle rerouting, \citep{bayraktar2023traffic} etc. Besides, the collected traffic data for each UAV detected location (link or node) become time series data considering the signal cycle, the UAV battery and the task rotation between multiple UAVs. Specifically, the under-determined system caused by the difference between limited UAV number and the network size poses a challenge for UAV-based traffic monitoring, that is how to know the network-wide traffic state with limited observation from UAVs. Mostly, such studies on network traffic monitoring mainly adopt probability distribution modeling for quantification of traffic state estimation error, or utilize neural networks to model the mapping between limited observation and joint feature of all network elements \citep{garcia2023trajectory, englezou2024enhancing, xiong2025multi}. Combined with multi-UAV routing or location optimization, a multi-objective model is proposed to minimize both the UAV operation cost and the traffic state estimation error \citep{salama2022collaborative, meng2023multi, theocharides2024real, bai2025dynamic, meng2025multi}. 

To sum up, UAV-based traffic monitoring is an emerging research branch closely related to the development of image processing and UAV scheduling, whose modeling idea and methodology are bottom-up and problem-oriented \citep{zandieh2023integrated, lin2025advanced}. Existing studies under signalized network scenarios focus more on the temporal correlation of traffic state than the spatial correlation among links and nodes in the network, which makes the traffic state of the network scenario just a combination of link-based or node-based observation or estimation, while neglecting the additional information carried by the topological relation between links and nodes. Moreover, with the full sample observation of UAV detection, detailed traffic dynamics like queue length, arrival rate, etc., are rarely studied.

\section{Preliminaries}\label{sec:research problem}

\subsection{Assumptions and definitions}

According to the signalized network scenario in Fig. \ref{fig: problem}, definitions on related terminology for network traffic measurement are given below.
\begin{itemize}
    \item \textbf{Network}. A network is a roadway system consisting of different levels of roads that connect to each other in a meshed topology. A tuple of graph $\mathcal{G} (\mathcal{I}, \mathcal{L})$ is used for notation, representing an intersection set $\mathcal{I}$ (indexed by $i$) and a link set $\mathcal{L}$ (indexed by $l$) respectively.
    \item \textbf{OD}. For a vehicle traveling within the network, its route has a traffic generation point (origin) and a traffic attraction point (destination), which are known as the OD and denoted by the corresponding links, $\mathcal{OD}={(\mathcal{L_O},\mathcal{L_D} )}$.
    \item \textbf{Path}. Topologically, a path refers to the sequential link set (or intersection set with OD links) for a vehicle to travel from one origin to one destination within the network, and the path set is denoted by $\mathcal{P}$ (indexed by $p$). From the perspective of signal control, a path can also be described by a sequential movement set regarding each intersection it passes. 
    \item \textbf{Movement}. For each intersection, the movements of left-turn, straight, and right-turn are generally considered for each connecting link; here, a set $\mathcal{M}_i$ is denoted to include all the movements of intersection $i \in \mathcal{I}$, indexed by $m$. 
    \item \textbf{Connectivity matrix}. For any two intersections $i,i'$ in network $\mathcal{G} (\mathcal{I}, \mathcal{L})$, if there is a link directly connecting two intersections, then the two intersections can be regarded as connected. A connectivity matrix $\bm{C} \in \mathbb R^{}$ is defined here for denotation, whose element $c_{i,i'}$ is given by Eq. \eqref{eq: connectivity}.
    \begin{align}
    c_{i,i'} = \begin{cases}
        1 \quad \text{if} \quad \exists \quad \text{link } l \text{ starts at $i$ and ends at $i'$} \\
        0 \quad \text{otherwise},
    \end{cases} \label{eq: connectivity}   
    \end{align}
\end{itemize}
For the following modeling part, some assumptions are made.
\begin{enumerate}
    \item The number of available UAVs $N_{uav}$ is fixed and they can all be used for traffic monitoring.  
    \item The UAV location to be optimized in this study refers to the center of any intersection, which the UAV is assigned to hover over. For such an intersection detected by UAV, the inflows and outflows of all the entries and exits are assumed observable. It also means that for the neighboring intersections, the upstream arrival flow of the connecting link is observable.
    \item For neighboring intersections both monitored by UAVs, vehicles passing the connected link are assumed to be matched through image processing \citep{coifman2006roadway}.
    \item If a certain movement $m$ is observable by both CVs and UAVs, CVs can be matched to the corresponding object in UAV detection. Similarly, if a certain movement $m$ is observable by both CVs and loop detectors, CVs can be matched to the corresponding passing timestamps in loop detection.
    \item The OD pair set $\mathcal{OD}$ and path set $\mathcal{P}$ are known \citep{behara2020novel}.
    \item The hovering height of UAVs at different intersections are assumed the same, thus the detection scale of UAVs at different intersections varies with the intersection topology, which affects the observability of traffic arrival and queue length.
    \item Path preference can be reflected through the historical data of CVs.
\end{enumerate}

\begin{figure}[!htbp]
\centering
\includegraphics[width=4.8 in]{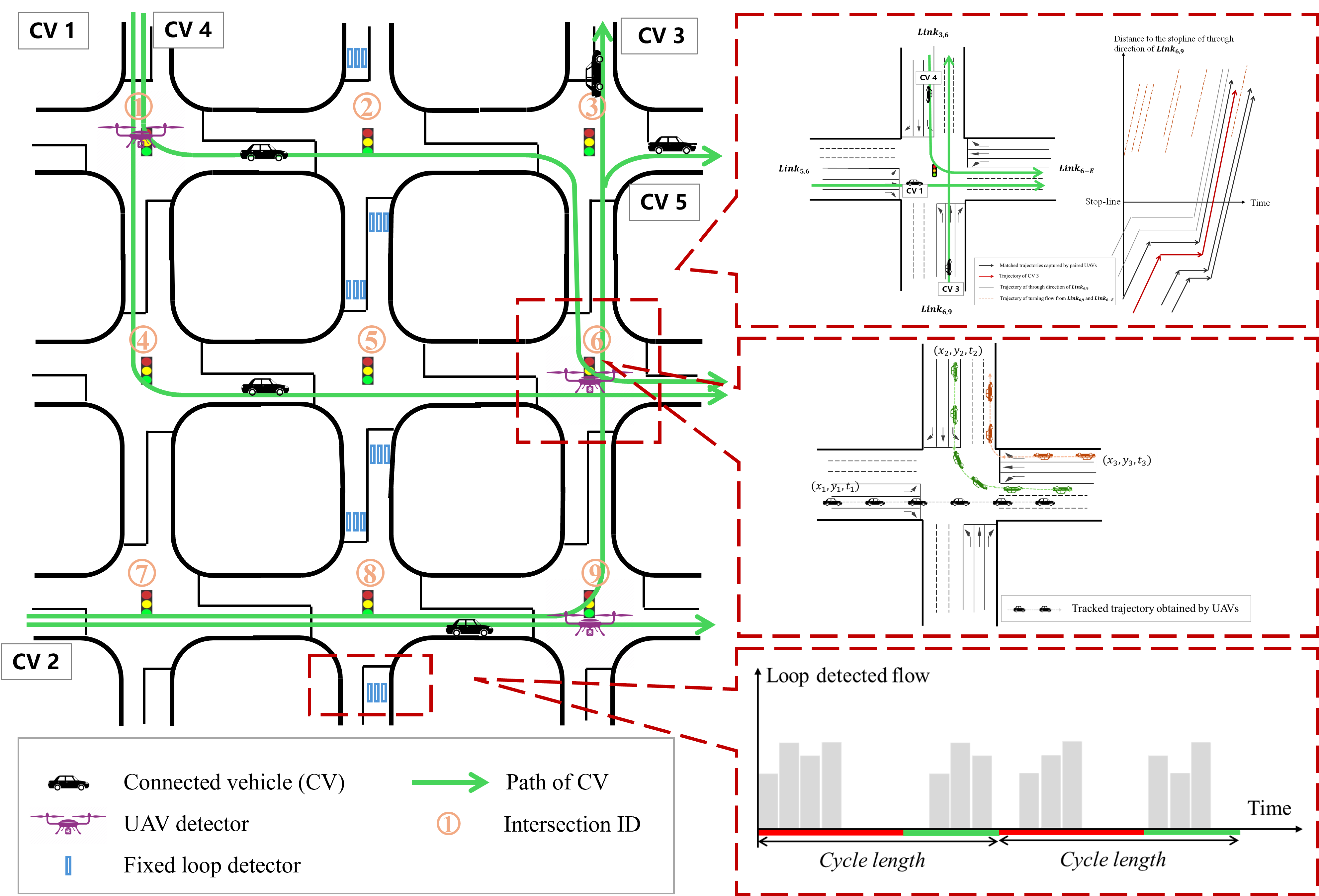}
\caption{\rmfamily\small A sketch of the research scenario.}
\label{fig: problem}
\end{figure}

\subsection{Problem statement}

Considering the multi-source detection environment, including both low-altitude and ground sensors shown in Fig. \ref{fig: problem}, network-wide traffic state can be partially observed at the network level and intersection level as stated below. 
\begin{itemize}
    \item From the network level, \textbf{an observed sub-path} refers to the sequence of movements of paired intersections with UAV sensors that a certain path $k$ passes through. Similarly, an observed global path refers to the sequence of movements of a complete CV trajectory along a certain path $k$. Such information facilitates the reconstruction of individual path choice, which means the actual path taken by a vehicle from a candidate path set $\mathcal{K}_{o}$ given the multi-source observation $o$.

    \item From the intersection level, \textbf{an observed traffic arrival} refers to the time-varying arrival rate curve of a certain movement $m$ at intersections detected by UAVs, also known as the arrival profile. As for links installed with lane-based loop detectors located about 20 to 30 meters upstream of the stopline, the passing timestamps obtained by event-actuated detection also reflect the distribution of traffic arrival among different movements. Compared with such full sample detection, sampled CVs, especially the queued ones when approaching the intersection, provide the partial arrival flow joining the queue before them or between them in the same cycle, which can be used as an inference of the traffic movement demand at the cycle level. 
    
   \item From the intersection level, \textbf{an observed queue length} refers to the back of the queue of each cycle of a certain movement $m$ at intersections detected by UAVs, while the position of queued CVs also offers the lower bound of the cycle-based queue length. Though the long queue problem of the fixed loop detectors can reflect the queue length to some extent, the case of no traffic arrival also shows no detection, similar to the case where the detector is occupied by a long queue. Thus, the observation of loop detectors should be integrated with other detection means for deterministic traffic index detection.
\end{itemize}

Thus, the task of network-wide traffic state monitoring can be concretized into three sub-tasks: individual vehicle path reconstruction, movement demand estimation, and queue length estimation. 
Compared with existing ground sensors like loop detectors or CVs, UAVs capture both the individual path choice and full vehicle trajectories at intersections, which features the consistency between microscopic travel behavior and macroscopic traffic flow parameters. Such an advantage thus facilitates the possibility of traffic monitoring that links network-wide traffic dynamics and individual path choice through UAV detection, which differs from existing studies that mainly focus on either movement-level or intersection-level traffic monitoring only.

Let $\mathcal{I}_{uav}$ be the subset of intersections deployed with UAVs, $\mathcal{L}_{fix}$ be the subset of links with fixed loop detectors installed, $\mathcal{T}_{cv}$ be the set of sampled CVs. Network-wide traffic monitoring is thus partially known from such multi-source observation $O$. Such imperfect knowledge of network-wide traffic monitoring is due to the limited coverage of detectors, which is inevitable in practice and is bound to cause uncertainty in traffic state estimation. Here, the uncertainty refers to the reliability of a random variable, which can be quantified by the feasible zone size given available detection observation. However, most current traffic state estimation methods focus on how to decrease the estimation error from the methodological level without considering the uncertainty shaped from the perspective of data input. In view of the dynamic detection property of CVs and UAVs, such uncertainty of traffic state actually can be obtained as a means for data quality evaluation and even decreased through detector configuration optimization, which helps reach more accurate estimation regardless of which method is used for estimation.

Therefore, the key challenge of this study lies in establishing the optimization model for UAV configuration driven by uncertainty of traffic monitoring, considering the consistency between intersection level and network level, which are handled through the three above-mentioned sub-tasks stated below.

\begin{itemize}
    \item From the network level, \textbf{vehicle path reconstruction task} can be modeled as a path selection inference problem. The uncertainty of vehicle path reconstruction facilitates the calculation of the path choice set, ultimately identifying the path with the highest probability as the true travel path.

    \item From the intersection level, \textbf{traffic movement demand and queue length estimation task} can be modeled as a fragmented curve fitting problem. The difference is that traffic movement demand estimation is to reconstruct the cumulative arrival and departure curves of a certain movement of the link, while queue length estimation is to reconstruct the space-time curves of individual vehicles, or shockwave curves generated by the status change of vehicles. Therefore, the uncertainty of traffic movement demand estimation facilitates the definition of the arrival profile uncertainty zone, which outlines the unknown part of the cycle-based traffic arrival curve so that the gap between traffic movement demand estimation and ground truth can be minimized. Similarly, the definition of back of queue (BoQ) uncertainty area is derived from the uncertainty of queue length estimation, which describes the possible location of where the queue propagates.
\end{itemize}

Using a uniform measurement of uncertainty, the UAV location problem can be formulated as finding an optimal selection of $\mathcal{I}_{uav}$ to minimize the above uncertainty of traffic state monitoring within a time period. Such an approach not only provides a paradigm of integrating low-altitude UAV detection and ground sensor detection for unified uncertainty formulation but also demonstrates generality across different networks. 

\section{Methodology}\label{sec:methodology}

\subsection{Location optimization model}
Given limited UAVs with prior sample CVs and fixed loop detectors as auxiliaries, UAV location optimization aims to minimize the uncertainty $Z$ of network-wide traffic state monitoring, as given by the multi-objective programming model:
\begin{align}
    \min_{I_{uav}} \quad & Z=w_1 F_{path} + w_2 F_{arrival} + w_3 F_{queue}. \label{eq: obj}\\
    \text{s.t.} \quad & u_i = \begin{cases}
        1 \quad \text{if} \quad i \in \mathcal{I}_{uav} \\
        0 \quad \text{otherwise},
    \end{cases} \label{eq:UAV location} \\
    & \sum_{i \in \mathcal{I}} u_i \leq N_{uav}.
\end{align}

As the paths can be regarded as a movement sequence and the intersection-level indicators are movement-based, the observability of the intersection should be transformed into the observability of movements. 

Here, a binary variable $y_{i,m}$ is defined to denote whether movement $m$ at intersection $i$ is observed by UAVs, as given by Eq.\eqref{eq: movement observability}. 
\begin{align}
    & y_{i,m} = \begin{cases}
        1 \quad \text{if} \quad \text{movement $m$ of intersection $i$ is observable by UAV} \\
        0 \quad \text{otherwise},
    \end{cases} \label{eq: movement observability}\\
    & y_{i,m} \geq u_{i} \label{eq: ego observability}\\
    & y_{i,m} \geq u_{i'}  \label{eq: upstream observability}\\ 
    & y_{i,m} \leq u_{i} + u_{i'} \label{eq: simultaneous observability} \\
    & \quad m \in \mathcal{M}_{i'-i},\quad  i'\in  \mathcal{I}_{adj}(i)
      \label{eq: interconnecting constraint}    
\end{align}
where  $\mathcal{I}_{adj}(i)$ represents the neighborhood intersection set of intersection $i$, thus $\forall  i'\in \mathcal{I}_{adj}(i), \ c_{i',i} = 1$.  $\mathcal{M}_{i'-i}$ denotes the movement set of the link starts at $i'$ and ends at $i$. Eqs. \eqref{eq: ego observability} - \eqref{eq: interconnecting constraint} declare that a certain movement can be regarded as observable by a UAV if a UAV is hovering at the intersection it belongs to or its upstream intersection. According to assumption 3, when both $u_i$ and $u_{i'}$ equal 1, then movement $m$ can be regarded as completely observable by the UAV.

The three sub-objectives correspond to the uncertainty of path reconstruction $F_{path}$, traffic movement demand $F_{arrival}$, and queue length $F_{queue}$ as mentioned in Section \ref{sec:research problem},  while $w_1$, $w_2$, and $w_3$ are weighting coefficients, respectively. For the following subsections, the formulations of $F_{path}$, $F_{queue}$, and $F_{arrival}$ are introduced in detail. 

\subsection{Uncertainty of path reconstruction}\label{sec:path reconstruction uncertainty}
In this section, an indicator, path choice uncertainty $U_{path}$ is proposed to quantify the uncertainty of path reconstruction considering the multi-source detection observation. 
It is noted that fixed loop detector data are not used in this subsection, as they cannot provide individual path information.
Considering $k$ as a random variable denoting the true path taken for an individual observed vehicle, and the global path set of the network is denoted as $K$. Without any detection data, the path choice set of a vehicle is $n$, which is the size of $K$.
For a user who actually chooses path $k$ for his or her trip, only a probability of $\frac{1}{n}$ can be sure that an unobserved vehicle chooses path $k$, thus its path choice uncertainty can be reflected by the size of the feasible path set $n$ to some degree. An uncertainty indicator is defined here to denote the probability of wrong path reconstruction, which is given as 
\begin{align}
    U_{path} = 1-\frac{1}{n}.
\end{align}

Let the random variable $O$ represent the observed sub-path set given a certain UAV location scheme. For sub-path $o$ captured by paired UAVs, $K_o$ represents the set of paths that pass through the observed sub-path $o$, with its size denoted as $n_o$. Similarly, only a probability of $\frac{1}{n_o}$ can be sure that an observed vehicle passing sub-path $o$ chooses path $k$, thus its path choice uncertainty can be formulated as
\begin{align}
    U_{path} = 1-\frac{1}{n_o}.
\end{align}

If CVs are captured for a UAV-detected sub-path $o$, the sampled CV flow can be used to model the path preference of $K_{o,cv}$. The path choice uncertainty of a vehicle passing sub-path $o$ is given by 
\begin{align}
    U_{path} = 1-p_{k,o} = 1-\frac{f_k}{f_o},
\end{align}
where $p_{k,o}$ denotes the path observation probability of path $k$ under observation $o$, $f_k$ denotes the CV path flow of path $k$, and $f_o$ denotes the CV flow of sub-path $o$. 

As for paths only observed by CVs, denoted as $K_{cv}$, the path choice uncertainty of a vehicle choosing one path from $K_{cv}$ is given by 
\begin{align}
    U_{path} = 1-\frac{f_k}{f_{cv}},
\end{align}
where $f_{cv}$ is the number of CVs excluding those passing any sub-path $o \in O$.

Thus, the path choice uncertainty of the whole network given different observation cases, given UAV and sampled CV data is,
\begin{align}
    F_{path} = \sum_{k \in \mathcal{K}} U_{path,k} \label{eq: path uncertainty}  
\end{align}
\begin{align}
    = \sum_{k \in \mathcal{K}_{o,cv}} [\frac{f_o}{Q_o}*(1-\frac{f_k}{f_o})] +  \sum_{k \in \mathcal{K}_o-\mathcal{K}_{o,cv}} (1-\frac{1}{n_o}) + \sum_{k \in \mathcal{K}_{cv}-\mathcal{K}_{o,cv}}(1-\frac{f_k}{f_{cv}}) + \sum_{k \in \mathcal{K}-\mathcal{K}_o-\mathcal{K}_{cv}} (1-\frac{1}{n_{non}}) \label{eq: sum of path uncertainty }
\end{align}
Here, the subscript $k$ is added in $U_{path,k}$ to denote the uncertainty of path $k$. The path set for path choice uncertainty calculation can be divided into four types, 1) the paths observed by both CVs and UAVs, $\mathcal{K}_{o,cv}$, 2) the paths only observed by UAVs, $\mathcal{K}_{o} - \mathcal{K}_{o,cv}$, 3) the paths only observed by CVs, $\mathcal{K}_{cv} - \mathcal{K}_{o,cv}$, 4) the paths neither observed by UAVs nor CVs, $\mathcal{K}-\mathcal{K}_{o}-\mathcal{K}_{cv}$, whose size is denoted by $n_{non}$. Based on the definition of path and movement in Section \ref{sec:research problem}, these four path sets can be obtained according to the observability of movements, as given by Eq. \eqref{eq: path observability 1} - \eqref{eq: two-fold observability}.

\begin{align}
    & card(\mathcal{K}_{o}) = \sum_{k} y_{k}^{o} \label{eq: path observability 1}\\
    & y_{k}^{o} \geq y_{i,m},\quad  \forall  m \in \mathcal{M}_k \label{eq: path observability 2}\\
    & y_{k}^{o} \leq \sum_{m \in \mathcal{M}_k} y_{i,m} \label{eq: path observability 3}\\
    & card(\mathcal{K}_{cv}) = \sum_{k} y_{k}^{cv} \label{eq: CV observability}\\
    & y_{k}^{cv} = \begin{cases}
        1 \quad \text{if any CV travels along path $k$} \\
        0 \quad \text{otherwise},
    \end{cases} \\
    & \mathcal{K}_{o,cv} = \mathcal{K}_{o} \cap \mathcal{K}_{cv} \label{eq: two-fold observability}
\end{align}

Here, $card( )$ denotes the size of a set. Eq. \eqref{eq: path observability 1}, \eqref{eq: CV observability}, and \eqref{eq: two-fold observability} describe the deduction of the path sets observable by UAV data, CV data and both data sources, respectively. $y_{k}^{o}, y_{k}^{cv}$ are 0-1 binary variables indicating whether path $k$ is observed by UAV or CV data respectively. As long as one movement of path $k$ is captured by U,AV then path $k$ can be regarded as observable, so constraints of Eq. \eqref{eq: path observability 2} - \eqref{eq: path observability 3} together determine the UAV observability of path $k$. $\mathcal{M}_{k}$ denotes the movement set path $k$ passes.

Fthe or path set $\mathcal{K}_{o,cv}$, a normalization step is used to involve the sub-path observation probability, $p(o)= \frac{f_o}{Q_o}$, to represent the relative magnitude of different sub-paths in terms of traffic demand. $Q_o$ denotes the aggregated flow of all observed sub-paths. 

Easily, we can find that $u_{path}< 1$ always exists.  A smaller $F_{path}$ indicates a smaller uncertainty of determining the true path for a vehicle of any path, which only relies on the UAV location rather than what kind of travel behavior assumption is adopted.

\subsection{Uncertainty of traffic movement demand}\label{sec:traffic demand uncertainty}
Vehicle arrivals are generally quantified as the arrival profile. In this section, the feasible domain size of arrival profile, $U_{arrival}$ is proposed to quantify the uncertainty of traffic movement demand given the traffic arrival between consecutive queued CVs, trajectories from UAV detection, and aggregated flow from fixed detectors. 

Arrival rate indicates the number of arrival vehicles per second in a cycle; thus, the global feasible area of arrival rate $S_{global,a}$ should be a rectangle in the arrival rate-time plane:
\begin{align}
    S_{global,\lambda} = \lambda_u C,
\end{align}
where $\lambda_u$ is the maximum arrival rate and $C$ is the cycle length.

The information provided by different data sources contributes to the shape of the arrival profile uncertainty zone to different degrees.
\begin{itemize}
    \item For intersections detected by UAVs, the trajectories of vehicles approaching the stopline, passing the intersection, and leaving for the exit can be observed partially due to the detection scale of UAVs. Thus, the exact movement-based space-time diagram can be obtained for the arrival profile reconstruction. For intersections whose neighboring intersections are detected by UAVs, the trajectories of the movements of the connecting link are partially observed by UAVs.
    \item For movements not observable by UAVs, only ground sensors can provide traffic flow information for the estimation of traffic movement demand. The queuing process of CVs at intersections can provide partial vehicle arrival information during the cycle. For example, by using the queuing positions and expected arrival time of two consecutive queued CVs, we can easily estimate the arrivals of non-CVs in between. As for links installed with fixed loop detectors, the arrival profiles of their movements are constrained by the aggregated flows of the lane-based loop detectors.
\end{itemize}

\subsubsection{Traffic movement demand uncertainty modeling given UAV deployment}
In the following text, determination of the feasible domain size of traffic movement demand will be discussed at the movement level regarding four types of observability described by $u_i$ and $u_{i'}$ based on the definition of $y_{i,m}$ in Eqs. \eqref{eq: ego observability} - \eqref{eq: interconnecting constraint}: (1) $u_i=1, u_{i'}=1$, (2) $u_i=1, u_{i'}=0$, (3) $u_i=0, u_{i'}=1$ and (4) $u_i=0, u_{i'}=0$. The type index will be added as the superscript of $U_{arrival}$ to denote the exact observability type, and so is the case for the queue uncertainty calculation in Section \ref{sec:queue length uncertainty}.
\begin{figure}[!htbp]
\centering
\includegraphics[width=0.8\textwidth]{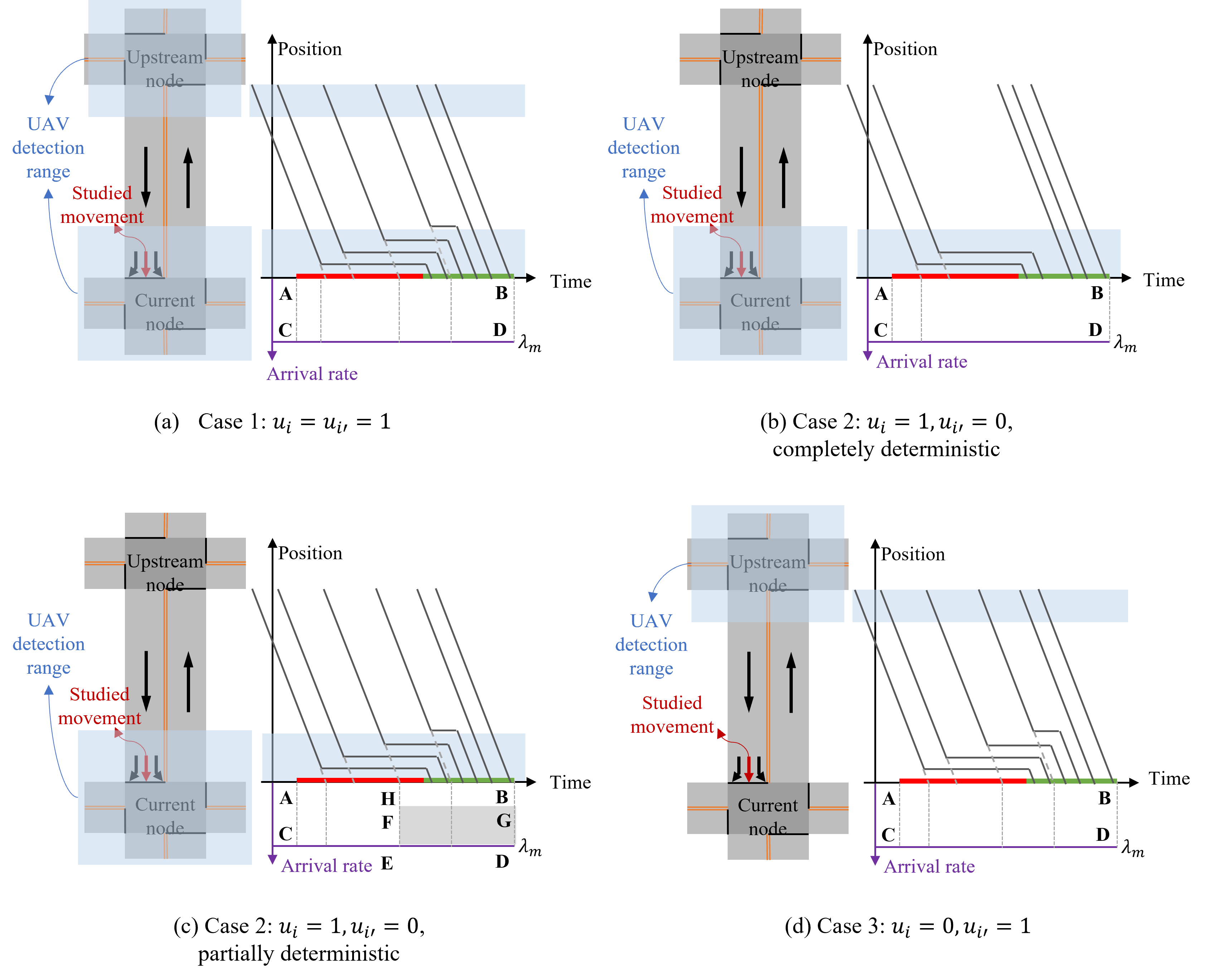}
\caption{\rmfamily\small Feasible area of vehicle arrivals for Case 1-3.}
\label{fig: UAV-based uncertainty of arrival}
\end{figure}

\textbf{Case 1: $u_i=1, u_{i'}=1$.}
Based on assumption 3, for two neighboring intersections detected by UAVs, the vehicle arrivals of all movements on the connecting link are fully deterministic, as shown in Fig. \ref{fig: UAV-based uncertainty of arrival}(a), the traffic movement demand uncertainty is $S_{u,\lambda} = 0$.

\textbf{Case 2: $u_i=1, u_{i'}=0$.}
For an intersection with a UAV hovering above, if the back of the queue is within the detection scale, then the vehicle arrivals within the cycle are fully known, as shown in Fig. \ref{fig: UAV-based uncertainty of arrival}(b). For such a case, the arrival rate uncertainty is $S_{u,\lambda} = 0$.

If the back of queue propagates beyond the detection scale, then the arrival profile during the cycle is partially determined. As shown in Fig. \ref{fig: UAV-based uncertainty of arrival}(c), the arrival rates during $T_{AH}$, are known, while arrival rates during $T_{HB}$ have a lower bound $\lambda_{HB,lb}=n_{HB}/T_{HB}$. Then, the feasible domain size of the arrival rate is 
    \begin{align}
        S_{u,\lambda} =(\lambda_u - \lambda_{HB,lb} )(t_G-t_F)
    \end{align}   

It is noted that when CVs and loops are also available, the arrival profiles of different movements can be further distinguished. 

\textbf{Case 3: $u_i=0, u_{i'}=1$.}
For an intersection with UAV observing at the neighboring intersection, the inflow of the movements on the connecting link can be obtained, though different movements share the same arrival profile by average, as shown in Fig. \ref{fig: UAV-based uncertainty of arrival}(d). It is noted that when fused with available CVs and loops, the arrival profiles of different movements can be further distinguished using the CV flow and aggregated loop flow of different movements, which can realize less uncertainty than simply using UAV data. 

\textbf{Case 4: $u_i=0, u_{i'}=0$.}
In case no UAVs capture the traffic flow information of a certain movement, the calculation of the feasible domain size depends on the ground sensor data, which is discussed according to the data source as below. And cases where both CV and loop data are available can be regarded as the combination of the single-source cases, and the uncertain areas of arrival rate can be correspondingly calculated.

\begin{figure}[!htbp]
\centering
\includegraphics[width=0.95\textwidth]{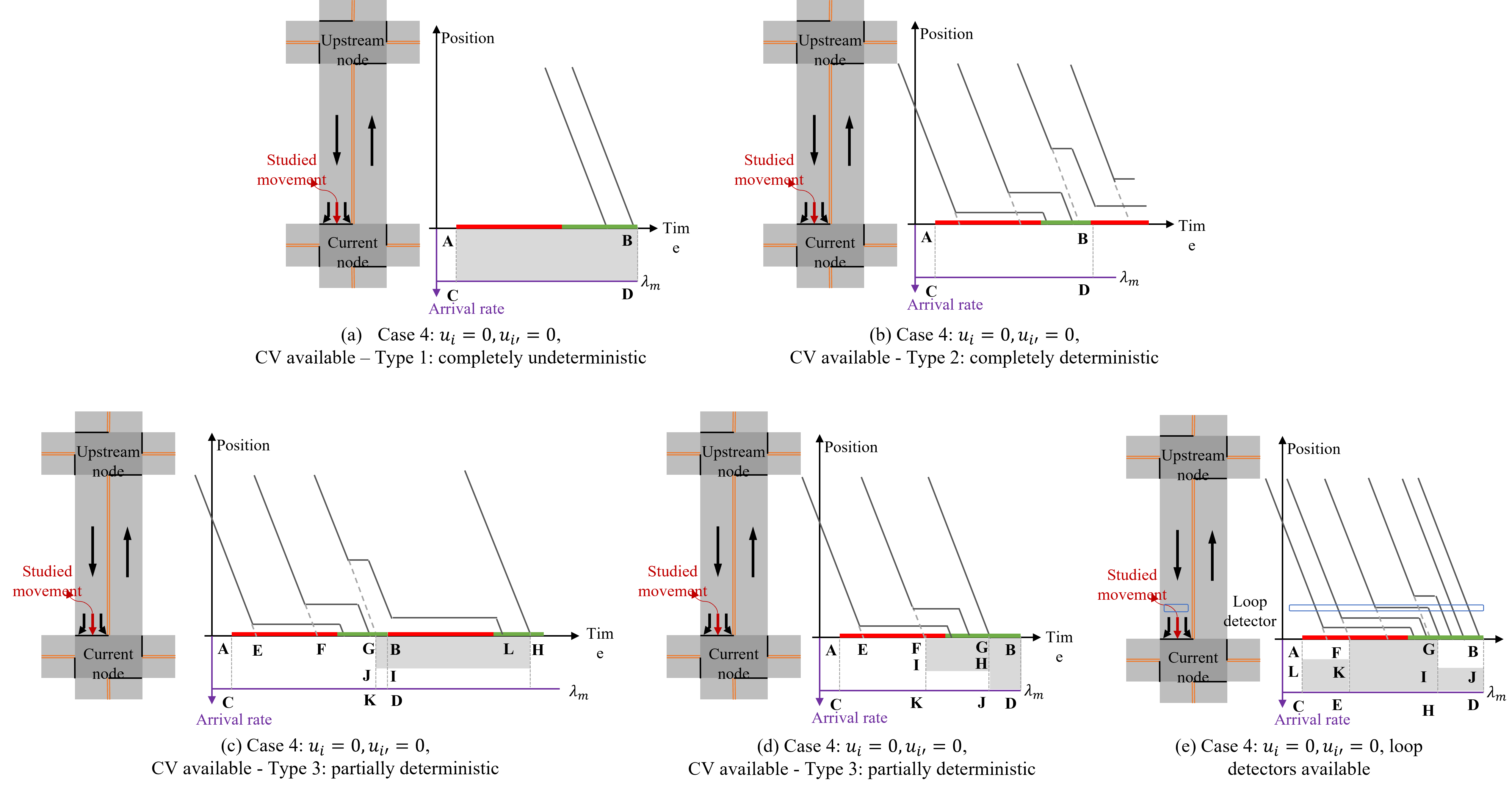}
\caption{\rmfamily\small Feasible area of vehicle arrivals for Case 4.}
\label{fig: uncertainty of arrival}
\end{figure}

For both undersaturated and oversaturated cycles, the feasible area of arrival rate can be cut based on CV trajectories, as the arrival rate between any combination of queued CVs (either queued or twice queued) is deterministic. Therefore, the arrivals at the beginning of the cycle are determined in the presence of CVs who queued for the first time in this cycle, while the arrival after the queue clearance during the green phase is determined by the presence of non-queued CVs or queued CVs left as residual queue till the next cycle. In summary, there are three types of feasible areas of arrival rates given CV observations during the cycle, as shown in Fig. \ref{fig: uncertainty of arrival}. It is noted that in Fig. \ref{fig: uncertainty of arrival} the vehicle trajectories are CV trajectories, which is different from the trajectories of the full sample of traffic flow in Fig. \ref{fig: UAV-based uncertainty of arrival}.

\begin{itemize}
    \item \textit{Type 1: Completely undetermined.} If no CVs are observed during the cycle, then the arrival rate of this cycle is completely indeterminate, as shown in Fig. \ref{fig: uncertainty of arrival} (a). The feasible area is 
    \begin{align}\label{eq: fully uncertain}
        S_{u,\lambda} = S_{global,\lambda}. 
    \end{align}
    \item \textit{Type 2: Completely determined.} If a cycle has twice-queued CVs (who queued for the first time in this cycle) and the next cycle has queued CVs, as shown in Fig. \ref{fig: uncertainty of arrival} (b), then the vehicle arrivals between all queued CVs are known, so the uncertainty is zero. 
    \begin{align}\label{eq: fully certain}
        S_{u,\lambda} = 0.
    \end{align}
    \item \textit{Type 3: Partially determined.} Except for Type 1 and 2 cycles, all the rest of the cycles belong to type 3, whose arrival profile during the cycle is partially determined. Generally, there are two subtypes.
    
    (i) As shown in Fig. \ref{fig: uncertainty of arrival} (c), a twice-queued CV is captured in the cycle, and only non-queued CVs are captured in the next cycle.     
    In this case, arrival rates during $T_{AE}$, $T_{EF}$, and $T_{FG}$ are known, while arrival rates during $T_{GH}$ have an upper bound $\lambda_{GH,ub}=T_{HL}/(T_{GH} * h_s )$. Then, the feasible domain size of the arrival profile is 
    \begin{align}
        S_{u,\lambda} = \lambda_{GH,ub} (C-t_G).
    \end{align}
    (ii) As shown in Fig. \ref{fig: uncertainty of arrival} (d), both queued and non-queued CVs are captured in the cycle. In this case, arrival rates $T_{AE}$ and $T_{EF}$ are known, arrival rates during $T_{FG}$ has an upper bound $\lambda_{FG,ub}=T_{GL}/(T_{FG} * h_s )$, and arrival rates during $T_{GB}$ is completely unknown, thus we have
    \begin{align}
        S_{u,\lambda} =\lambda_u (C-t_G )+\lambda_{FG,ub} (t_G-t_F )
    \end{align}
\end{itemize}

As the lane-based loop detectors are fixed, the arrival profiles of the movements of interest are partially determined, as a long queue problem will make them invalid for data collection. As shown in Fig. \ref{fig: uncertainty of arrival} (e), the vehicle count obtained during the interval $T_{AF}, T_{GB}$ when the detector is not occupied provides a lower bound for the arrival rate $\lambda_{AF,lb}=n_{AF}/T_{AF}, \lambda_{GB,lb}=n_{GB}/T_{GB}$, while the arrival rate within the period $T_{FG}$ when the detector is occupied is fully unknown. Thus, the feasible domain size of the arrival rate is 
    \begin{align}
        S_{u,\lambda} =\lambda_u*(t_G-t_F) + (\lambda_u - \lambda_{AF,lb})(t_F-t_A) + (\lambda_u - \lambda_{BG,lb})(t_B-t_G)
    \end{align}

\subsubsection{Traffic movement demand uncertainty formulation}

Based on the detection condition of a certain movement and the detection data within a cycle, the uncertainty index for arrival rate, i.e., $U_{arrival}$, can be given as below
\begin{align}
    U_{arrival} = S_{u,\lambda}/S_{global,\lambda}
\end{align}
We can find that $0 \leq U_{arrival} \leq 1$. Note that, $U_{arrival}=0$ does not mean that estimates by existing methods are error-free. 
A smaller $U_{arrival}$ indicates less uncertainty of the arrival profile, suggesting that any method is expected to achieve a more accurate arrival rate estimate based on the observed CV trajectories.

In summary, considering the above four observability types, the traffic movement demand uncertainty of the whole network is calculated as
\begin{equation}
\begin{split}
    F_{arrival} = \sum_{i \in \mathcal{I}} \sum_{m \in \mathcal{M}_i} \{y_{i,m} [u_{i}*u_{i'}\sum_{j \in \mathcal{J}_m} U_{arrival}^{1,i,m,j}+u_{i}*(1-u_{i}*u_{i'})\sum_{j \in \mathcal{J}_m} U_{arrival}^{2,i,m,j} + \\
    u_{i'}*(1-u_{i}*u_{i'})\sum_{j \in \mathcal{J}_m} U_{arrival}^{3,i,m,j}]+(1-y_{i,m})\sum_{j \in \mathcal{J}_m} U_{arrival}^{4,i,m,j} \} \label{eq: arrival uncertainty},
\end{split}
\end{equation}
where $\mathcal{M}_i$ denotes the movement set of intersection $i$, $\mathcal{J}_m$ denotes the set of all cycles of movement $m$ during the analysis period, $U_{arrival}^{*,i,m,j}$ denote the arrival rate uncertainty of cycle $j$ of movement $m$ at intersection $i$ for a certain observability type. 

\subsection{Uncertainty of queue length}\label{sec:queue length uncertainty}
Queue length are generally evaluated through the maximal length of BoQ within each cycle. In this section, the feasible domain size of BoQ, $U_{queue}$ is proposed to quantify the uncertainty of queue length considering the four observability types stated in Section \ref{sec:traffic demand uncertainty}. 

For a specific movement of a signalized intersection, the feasible domain size of the queue length within a signal cycle is dependent on the signal timing parameters. Given the maximum queue accumulation wave speed $w_a$, and the queue dissipation wave speed $w_d$, we can define a global feasible area $S_{global,q}$ of possible queue length, as the triangle area in the space-time diagram bounded by the stop-line, the curve denoting the queuing shockwave, and the dissipating shockwave: 
\begin{align}
    S_{global,q} = 0.5 w_a w_d R^2/(w_d-w_a)
\end{align}
where $R$ is the red time of the cycle.

\begin{figure}[!htbp]
\centering
\includegraphics[width=0.8\textwidth]{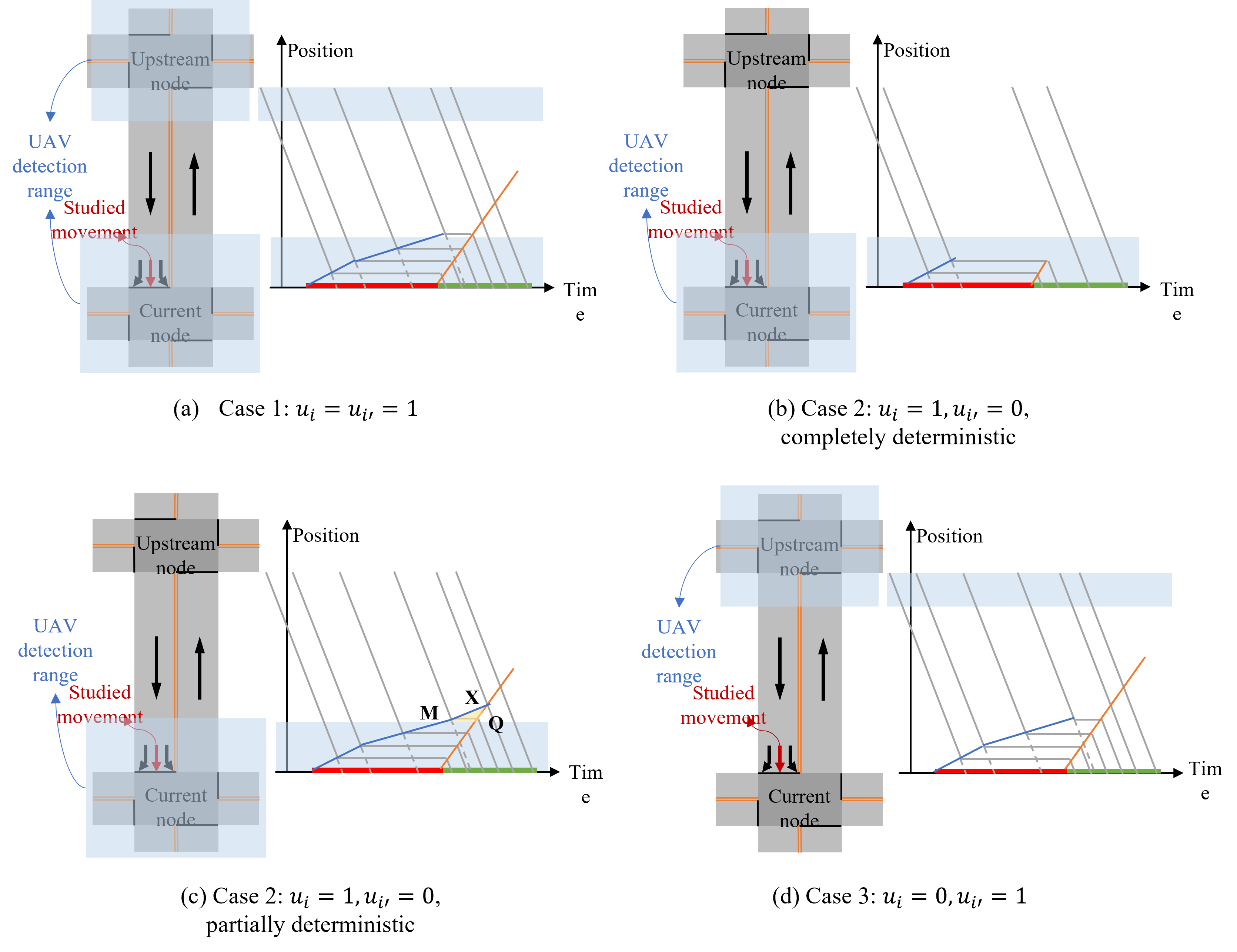}
\caption{\rmfamily\small Feasible area of queue length for Case 1-3.}
\label{fig: UAV-based uncertainty of queue}
\end{figure}

\subsubsection{Queue uncertainty modeling given UAV deployment}

\textbf{Case 1: $u_i=1, u_{i'}=1$.}
For two neighboring intersections detected by UAVs, the queue lengths of all movements on the connecting link are fully deterministic, as shown in Fig. \ref{fig: UAV-based uncertainty of queue}(a). For such cases, the queue uncertainty is $S_{u,q}=0$.

\textbf{Case 2: $u_i=1, u_{i'}=0$}
For an intersection with UAV hovering above, if the back of queue is within the detection scale, then the queue length within the cycle are fully known, as shown in Fig. \ref{fig: UAV-based uncertainty of queue}(b), $S_{u,q}=0$.  

When the queue propagates beyond the UAV detection scope, the feasible area of queue length is a triangle, as shown in Fig. \ref{fig: UAV-based uncertainty of queue}(c). Point $Q:(t_Q,d_Q )$ and Point $M:(t_M,d_M)$ are known based on the UAV detection bound and the farthest queued trajectory from the stop-line captured by UAV. Then the uncertainty area can be obtained using Eq. \eqref{eq: UAV partial queue uncertainty}. 
\begin{align}\label{eq: UAV partial queue uncertainty}
    S_{u,q}=0.5(t_Q-t_M )(d_X-d_Q ).
\end{align}

It is noted that with CVs and loops available, the queue length feasible zones of different movements can be further distinguished using the queued CVs and the position of loop detectors of different movements. 

\textbf{Case 3: $u_i=0, u_{i'}=1$}
For an intersection with UAV observing at the neighboring intersection, the inflow of the movements on the connecting link can be obtained, thus the queue length can be obtained, as shown in Fig. \ref{fig: UAV-based uncertainty of queue}(d). Similar to the uncertainty of arrival rate, the queue lengths of different movements can be further distinguished based on the calibrated turning ratio when fused with CVs and loop. 

\begin{figure}[!htbp]
\centering
\includegraphics[width=0.8\textwidth]{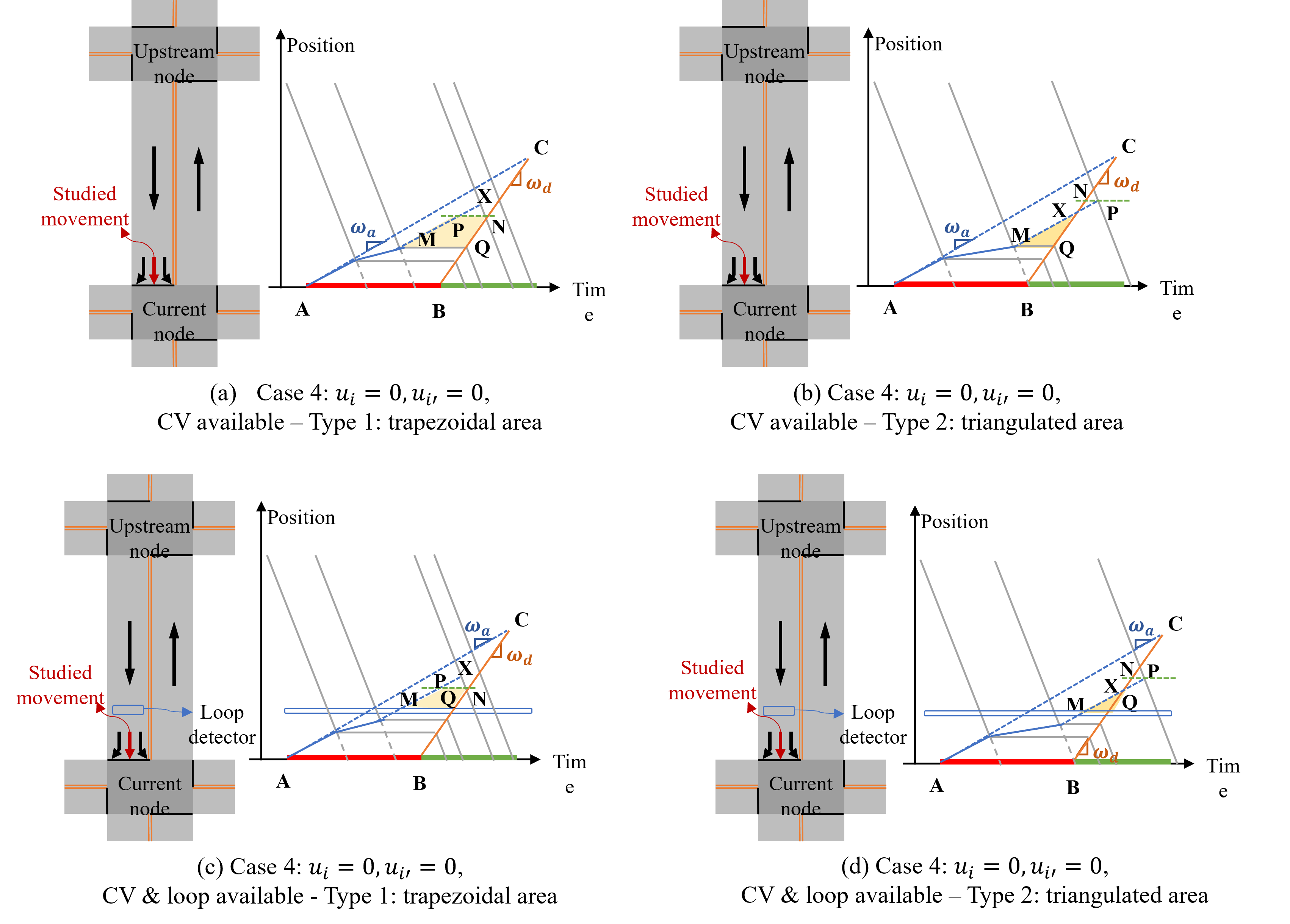}
\caption{\rmfamily\small Uncertain area of queue length for Case 4.}
\label{fig: uncertainty of queue}
\end{figure}

\textbf{Case 4: $u_i=0, u_{i'}=0$}
In case that no UAVs capture the traffic flow information of a certain movement, the calculation of feasible domain size depends on the ground sensor data, which is discussed according to the data source as below, shown in Fig. \ref{fig: uncertainty of queue}. It is noted that in Fig. \ref{fig: uncertainty of queue} the vehicle trajectories are CV trajectories, which is different from the trajectories of the full sample of traffic flow in Fig. \ref{fig: UAV-based uncertainty of queue}. And the case where both CV and loop data are available can be regarded as the combination of the single-source cases and the uncertain areas of queue length can be correspondingly calculated.

As the vehicle arrival of a movement may consist of several incoming traffic flows from different approaches of the upstream intersection, the actual queuing shockwave may be more often piecewise. Thus, whether the queue length can be deterministic depends on whether queued CVs of different platoons can be captured to provide the corresponding arrival information. 

As shown in Fig. \ref{fig: uncertainty of queue}(a) and (b), the feasible queue length is determined by two points: the joining queue point $M$ of the last queued CVs and the crossing point $N$ of the line segment $\overline{BC}$ and the first non-queued CV’s trajectory (no matter whether this CV succeeds in passing in the current cycle). The coordinates of point $M$ are denoted by $M:(t_M,d_M )$, while the coordinates of point $N$ are denoted by $N:(t_N,d_N)$, both of which can be easily solved based on simple plane geometry. Note that for the unsaturated and oversaturated traffic conditions, the differences in the model derivation lies in the process of obtaining $M$ and $N$. The subsequent steps are identical.

Given point $M$, we can generate a line parallel to $\overline{AC}$ and intersecting with $\overline{BC}$ at point $X:(t_X,d_X )$, and we can solve that
\begin{align}
    &t_X=(w_d R-w_a t_M+d_M)/(w_d-w_a ), \\
    &d_X=(w_d w_a (R-t_X )+d_X)/(w_d-w_a ).
\end{align}

Then, based on the relative positions of $X$ and $N$, the feasible area of queue length has two types of shapes:

\begin{itemize}
    \item \textit{Type 1: trapezoidal area.} As shown in Fig. \ref{fig: uncertainty of queue} (a), when $d_X>d_N$, the feasible area of queue length, i.e., $S_{u,q}$, is a trapezoid. Point $Q:(t_Q,d_Q )$ is also known based on the trajectory of the last queued CV. Point $P:(t_P,d_P)$ is solved as $t_P=t_M+(d_N-d_M)/w_a$ and $d_P=d_N$. Then we have 
    \begin{align}\label{eq: trapezoidal queue uncertainty}
        S_{u,q}=0.5 (t_Q-t_M+t_N-t_P )(d_N-d_Q )
    \end{align}
    \item \textit{Type 2: triangulated area.} As shown in Fig. \ref{fig: uncertainty of queue}(b), when $d_X \leq d_N$, $S_{u,q}$ is a triangle. We have 
    \begin{align}\label{eq: triangle queue uncertainty}
        S_{u,q}=0.5(t_Q-t_M )(d_X-d_Q ).
    \end{align}
\end{itemize}

As loop detectors only collect the passing timestamps rather than the trajectories of the approaching vehicles, whether the passing vehicles join the queue or not can not be determined, so does whether the interval without data is caused by the queue or actually no vehicles pass the loop. Thus, the uncertainty area of the queue given loop observation can only be determined when fused with CVs, as shown in Fig. \ref{fig: uncertainty of queue} (c) and (d). Only when both queued and non-queued CVs are captured within the cycle, the location of the loop detector can help to reduce the uncertainty area of queue length if it lies between the CV queuing position and the calibrated back of queue. The uncertainty area of queue length can be obtained through Eq. \eqref{eq: trapezoidal queue uncertainty} and Eq. \eqref{eq: triangle queue uncertainty}, where $d_Q, d_M$ are the locations of the loop detector.

\subsubsection{Queue length uncertainty formulation}

Based on the detection condition of a certain movement and the detection data within a cycle, the uncertainty index for queue length, i.e., $U_{queue}$, can be given as below:
\begin{align}
    U_{queue} = S_{u,q}/S_{global,q}.
\end{align}

The physical meaning of $U_{queue}$ is the uncertainty regarding the spatial-temporal plane, even before any estimation model is applied. 
We can find that $0 \leq U_{queue} \leq 1$, and this index is independent of any theoretical model for queue length estimation, which is a direct evaluation of the input detection data. A smaller $U_{queue}$ indicates a smaller possible area of queue length, suggesting that any method is expected to achieve a more accurate queue length estimate based on the available detection data.

Considering the above four observability types, the queue length uncertainty of the whole network is calculated as
\begin{equation}
\begin{split}
    F_{queue} = 
     \sum_{i \in \mathcal{I}} \sum_{m \in \mathcal{M}_i} \{y_{i,m} [u_{i}*u_{i'}\sum_{j \in \mathcal{J}_m} U_{queue}^{1,i,m,j}+u_{i}*(1-u_{i}*u_{i'})\sum_{j \in \mathcal{J}_m} U_{queue}^{2,i,m,j} + \\
    u_{i'}*(1-u_{i}*u_{i'})\sum_{j \in \mathcal{J}_m} U_{queue}^{3,i,m,j}]+(1-y_{i,m})\sum_{j \in \mathcal{J}_m} U_{queue}^{4,i,m,j} \}
    \label{eq: queue uncertainty}
\end{split}
\end{equation}
where $U_{queue}^{*,i,m,j}$ denote the queue length uncertainty of cycle $j$ of movement $m$ at intersection $i$ for a certain observability type. 

\subsection{Solution algorithm}

Based on the uncertainty formulation in Section \ref{sec:path reconstruction uncertainty} - Section \ref{sec:queue length uncertainty}, the UAV location optimization model in Eq. \eqref{eq: obj} can be organized as below.
\begin{align}
    \min_{I_{uav}} \quad & Z=w_1 \sum_{k \in \mathcal{K}} U_{path,k} + w_2 \sum_{i \in \mathcal{I}} \sum_{m \in \mathcal{M}_i} y_{i,m} \sum_{j \in \mathcal{J}_m} U_{arrival}^{i,m,j} + w_3 \sum_{i \in \mathcal{I}}  \sum_{m \in \mathcal{M}_i} y_{i,m} \sum_{j \in \mathcal{J}_m} U_{queue}^{i,m,j} \label{eq: unfolded obj}\\
    & \text{s.t. Eq.\eqref{eq:UAV location} - Eq. \eqref{eq: queue uncertainty}}.
\end{align}

From the above formulation, the location of UAV, which is denoted by the intersection set $\mathcal{I}_{uav}$, actually determines the observed path $\mathcal{K}_{o}$ and the observed movement, which is defined as $\mathcal{M}_{o}$. It is noted that the size of $\mathcal{M}_{o}$ is larger than the movements of $\mathcal{I}_{uav}$ due to the UAV detection of traffic flows merging into the exits of $\mathcal{I}_{uav}$. From Eq. \eqref{eq: path uncertainty}, \eqref{eq: arrival uncertainty} and \eqref{eq: queue uncertainty}, the model of Eq. \eqref{eq: unfolded obj} is a non-linear integer programming model, which is actually an NP-hard problem and hard to solve using exact algorithms within an acceptable computation time. Especially for the research scenario in this study, with the increase of the roadway network scale, the problem presents nonlinearity with a large number of local extremes, which calls for better diversity representation and evolution randomness.

Quantum genetic algorithm (QGA) \citep{Narayanan1996, malossini2008quantum} is adopted here to obtain the solution integrating the thoughts of metaheuristic algorithms and quantum computation. As compared with classical GA, QGA adopts quantum bits for population encoding and quantum rotation gates for population evolution, which overcomes the shortcomings of easily falling into local maxima, slow convergence spee,d and poor accuracy of search results.

As the objective of the proposed method is to minimize the network-wide traffic state estimation uncertainty, here $-Z$ is used as the fitness function. 
Under the context of the research problem, the classic QGA is improved with the integration of two operators: a) the deduplication operator and b) the neighborhood-based fine-tune operator. The deduplication operator aims to improve the solution variety by removing the same individuals in each population evolution based on quantum mutation, while the fine-tune operator aims to explore a better solution by searching the neighborhood of the best individual of each generation based on the matching UAV configuration of assumption 3. With these two added operators, premature convergence can be avoided and the population diversity can be expanded, thus increasing the opportunity of searching better solutions.
The pseudocode of the proposed improved QGA (IQGA) is shown in Algorithm \ref{alg:IQGA}.

\begin{algorithm}
\caption{Improved Quantum Genetic Algorithm (IQGA)}
\label{alg:IQGA}
\begin{algorithmic}[1]
\REQUIRE Network topology $(\mathcal{I},\mathcal{L})$, sample CV data, UAV fleet size $N_{uav}$, chromosome population size $P$, chromosome length $I$ which is also the size of the intersection set, number of evolution generations $T$, number of mutated location $R_m$, boundaries of rotation angle $\theta_{min}, \theta_{max}$
\ENSURE Optimal UAV location $\mathcal{I}_{uav}$
\STATE Set generation counter $t=0$
\STATE Generate initial population: $P(0) = \{\text{Chrom}_1^{(0)}, \text{Chrom}_2^{(0)}, \dots, \text{Chrom}_{P}^{(0)}\}$. A quantum chromosome with length $I$ can be represented by qubits in Eq. \eqref{eq: quantum chromosome}, where $[\alpha^{t}_{pi} , \beta^{t}_{pi}]$ denote the probability amplitudes of the quantum state and satisfy the normalization condition $|\alpha|^2+|\beta|^2 = 1$. A qubit can thus be  expressed as a linear superposition of $ \ket{0}$ and $\ket{1}$, $\ket{\varphi} = \alpha \ket{0} + \beta \ket{1}$. 
\begin{align} \label{eq: quantum chromosome}
    Chrom_p^{(0)} = \left[
\begin{vmatrix} \alpha^{t}_{p1} \\ \beta^{t}_{p1} \end{vmatrix}
\begin{vmatrix} \alpha^{t}_{p2} \\ \beta^{t}_{p2} \end{vmatrix}
\dots
\begin{vmatrix} \alpha^{t}_{pI} \\ \beta^{t}_{pI} \end{vmatrix}
\right]
\end{align}
\STATE Measure the observation value of population $p$ and obtain the status $O(t,p)={o^{t,p}_1, o^{t,p}_2, \dots, o^{t,p}_I}$, 
\STATE Decoded $O(t)$ to obtain the UAV location scheme, $\vec{u}(t,p) = \{u_1^{t,p}, u_2^{t,p}, \dots, u_i^{t,p}, \dots\}$, indicating whether a UAV is set to hover above intersection $i$ for detection. 
\STATE Calculate the fitness value $\text{fit}(\vec{u}(t,p)) = -Z(\vec{u}(t,p))$
\STATE $B(t) = \arg\max_{p \in P} \text{fit}(\vec{u}(t,p))$
\STATE Record optimal individual $\hat{Chrom}$ and its fitness value  $\vec{u}^{t,\text{best}} = \arg\max_{p \in P} \text{fit}(\vec{u}(t,p))$
\WHILE{$t < T$}
    \STATE Evolve $t$ to $t+1$ using quantum rotation gate: $t = t + 1$. A fitness-based dynamic rotation strategy is adopted to determine the quantum rotation angle, as given by Eq. \eqref{eq: rotation angle}.
    \begin{align}
    \Delta \sigma_{t} = \theta_{min} + (\theta_{max}-\theta_{min})*\frac{|fit^{cur}-fit^{best}|}{max \{ fit^{cur},fit^{best} \} }
    \label{eq: rotation angle}
    \end{align}  
    For generation $t+1$, the quantum rotation angle updates as $\sigma_{t+1} = \sigma_{t} + \Delta \sigma_{t}$, and the $i$ th bit of the chromosome is updated as follows:
    \begin{align}
    \begin{pmatrix} \alpha^{t+1}_{pi} \\ \beta^{t+1}_{pi} \end{pmatrix} = 
    \begin{pmatrix} cos(\sigma_{t+1}) \quad -sin(\sigma_{t+1}) \\ sin(\sigma_{t+1}) \quad cos(\sigma_{t+1}) \end{pmatrix}
    \begin{pmatrix} \alpha^{t}_{pi} \\ \beta^{t}_{pi} \end{pmatrix}
    \label{eq: rotation gate}
    \end{align} 
    Obtain new population $P(t+1)$. 
    \STATE Duplicate checking. For the whole population, if the chromosome of population $p$ is the same as that of population $p'$, employ quantum NOT gate to realize quantum mutation, as given by Eq. \eqref{eq: quantum not gate} which actually changes the probability amplitudes of a certain qubit to flip the optimization direction of chromosomes. The qubits to be mutated are determined by randomly choosing $R_m$ from the chromosome length $I$.
    \begin{align}
    \begin{bmatrix} \alpha'^{t}_{pi} \\ \beta'^{t}_{pi} \end{bmatrix} = 
    \begin{bmatrix} 0 \quad 1 \\ 1 \quad 0 \end{bmatrix}
    \begin{bmatrix} \alpha^{t}_{pi} \\ \beta^{t}_{pi} \end{bmatrix}
    \label{eq: quantum not gate}
    \end{align}
    \STATE \textbf{Repeat lines 3 to 6.}          
    \IF{ $\text{fit}(B(t)) < \text{fit}(\vec{u}^{t,\text{best}})$}
        \STATE $\vec{u}^{t,\text{best}} = B(t)$, and the optimal individual $\hat{Chrom}$ is updated as the one with the new updated fitness.
        \STATE Fine-tuning based on neighborhood searching. 
        Replace the UAV deployment of one intersection with the minimal traffic movement demand and queue length uncertainty with one neighboring intersection of the UAV-detected intersection with the maximal traffic movement demand and queue length uncertainty, and evaluate the fitness value of the alternative scheme $\vec{u'}(t,p)$, if $\text{fit}(\vec{u'}^{t,p}) < \text{fit}(\vec{u}^{t,p})$, update the population by replacing the encoded scheme with the original one for the next generation.
    \ENDIF
\ENDWHILE
\end{algorithmic}
\end{algorithm}

\section{Evaluation}

The proposed method is evaluated based on an empirical network in Qingdao, China. Experiments are conducted in three aspects: 1) Performance of the proposed UAV deployment optimization method is validated regarding the relationship between objectives of multi-resolution traffic state uncertainty and the size of available UAV fleet, as well as sensitivity analysis on the balance between different sub-objectives, optimal UAV location scheme, and the optimal UAV fleet size. 2) Two types of baselines of UAV location methods, flow coverage maximum (FCM), and single-objective uncertainty minimum (SOUM) are selected for horizontal comparison, while the proposed method is referred to as multi-objective uncertainty minimum (MOUM). Three downstream applications, BoQ estimation, arrival rate estimation and path flow estimation are chosen to compare the effectiveness of different UAV location methods for traffic state monitoring. 3) The superiority of the proposed IQGA is justified through comparison with the classic solution algorithm, QGA, with regards to both exploitation and exploration ability.

\subsection{Study site}
\begin{figure}[!htbp]
\centering
\includegraphics[width = 5.6 in]{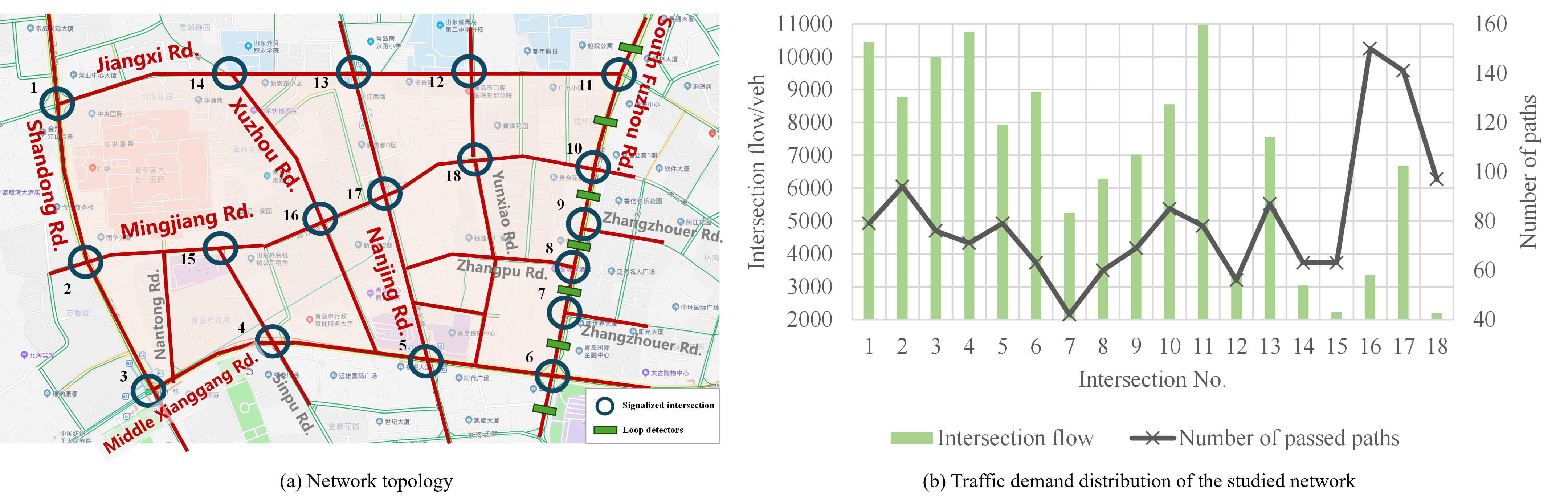}
\caption{\rmfamily\small Study site.}
\label{fig: shinan_network}
\end{figure}

As shown in Fig. \ref{fig: shinan_network}, an 18-intersection roadway network in Shinan District, Qingdao, Shandong Province, is selected as the study site. As marked out by blue circles with index in black, the study site includes 14 typical signalized crossroads and 4 T-type signalized intersections. All intersections are installed with license plate recognition (LPR) detectors, and lane-based fixed loop detectors are installed on the two-way mainline of South Fuzhou Road, as indicated by the green bar covering  6, 7, 8, 9, 10 and 11 in Fig. \ref{fig: shinan_network}(a).
Through VISSIM software, the simulation model of this Shinan network was built, and the empirical data from 7:00 – 9:00 on Mar. 1st, 2019, including LPR data, fixed loop detector data, floating car data (FCD), and signal timing data were used for calibration. Based on empirical detection data, there were in total 28 origins (destinations) and 311 paths, with the number of paths passing each intersection shown in Fig. \ref{fig: shinan_network}(b). It is noted that the axis label of \textit{Intersection flow} of the primary vertical axis refers to the sum of traffic flows passing a certain intersection in different movements, while the axis label of \textit{Number of passed paths} of the secondary vertical axis refers to the number of paths constituting any movements of a certain intersection.
    
The calibrated simulation model was run three times, with a simulation period of 2 hours each time. For the multi-source data input for UAV location optimization, a penetration rate of 0.1 was applied to extract the vehicle trajectories at an uploading frequency of 1s to obtain the sampled trajectories as the CVs, while vehicle counters were set at the two-way mainline of South Fuzhou Road to simulate the empirical loop detector configuration for aggregated flow data collection. As for the UAV detection, a detection scale of a 200m*200m box is assumed to be captured when a UAV is hovering over the center of an intersection. With the intersection scale varying, the number of arrival trajectories to the queue before the stopline during the red phase is different from the UAV detection perspective.

\subsection{UAV deployment evaluation}

\subsubsection{Traffic state uncertainty analysis}
Coded in Python, the method was run on a laptop with an AMD Ryzen 7 5800H CPU @ 3.20GHz, 16.0 GB of RAM, and one NVIDIA GeForce RTX 3060 GPU. For the IQGA solver, the population and generation were set as 20 and 200, respectively. Given the sampled CV data and loop detector as input, the proposed method was tested with $N_{uav}$ ranging from 1 to 17, with three parallel tests for each certain UAV number. Based on assumption 3, when the UAV coverage is 100\% (18 UAVs), the traffic state is fully determined, thus the uncertainty is 0 and not shown in the figure. 

\begin{figure}[!htbp]
\centering
\includegraphics[width=6.4 in]{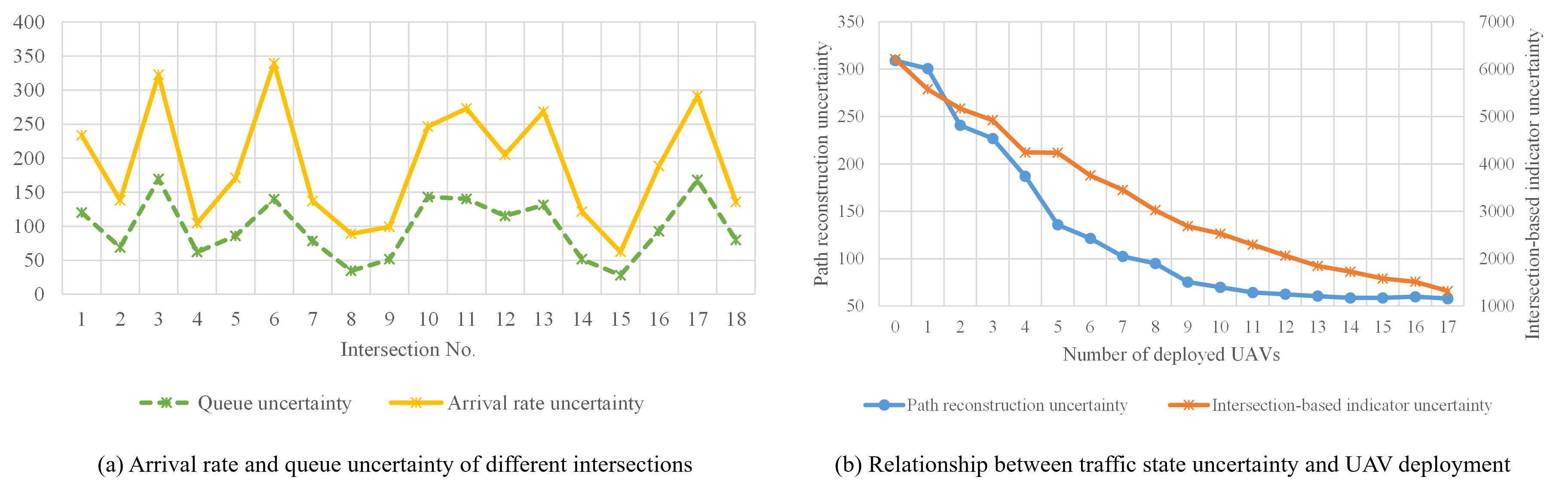}
\caption{\rmfamily\small Relationship between traffic state uncertainty and the number of deployed UAVs}
\label{fig: uncertainty trend}
\end{figure}

The values of three sub-objectives are averaged for three parallel tests to ensure generality. As the uncertainties of traffic movement demand and queue length are cycle-based, the values of sub-objectives $F_{arrival}$ and $F_{queue}$ are comparable and the distribution among different intersections is similar, as shown in Fig. \ref{fig: uncertainty trend}(a). Thus, the values of these two sub-objectives can be aggregated and denoted as intersection-based indicator uncertainty. It is noted that the queue uncertainty and arrival rate uncertainty in Fig. \ref{fig: uncertainty trend}(a) are calculated at the intersection level, stated as below.
\begin{itemize}
    \item \textbf{Uncertainty of queue length}.The feasible area to determine the exact back of the queue for all the movements of an intersection, which can be calculated by $\sum_{m \in \mathcal{M}_i} \sum_{j \in \mathcal{J}_m} U_{queue}^{i,m,j}$.
    \item \textbf{Uncertainty of arrival rate}.The feasible area to determine the arrival profiles for all the movements of an intersection, which can be calculated by $\sum_{m \in \mathcal{M}_i} \sum_{j \in \mathcal{J}_m} U_{arrival}^{i,m,j}$.
\end{itemize}

As the path reconstruction uncertainty is period-based, the weight coefficients are set as $w_1:w_2:w_3 = 26:1:1$ to ensure the magnitudes of individual vehicle path reconstruction, traffic movement demand estimation, and queue length estimation are nearly equal in network-wide traffic state uncertainty. Before that, a sensitivity analysis is applied to the weight coefficients as shown in Table \ref{tab: Sensitivity of weight coefficients}, which justifies the value setting of $w_1:w_2:w_3 = 26:1:1$ for further experiments, indicating the equal contribution of each sub-objective. It is noted that the sensitivity analysis is conducted for the case of $N_{uav}=9$ as the solution space of choosing 9 intersections from a total of 18 is the largest in this case study.

With the calibrated weight coefficients, Figure \ref{fig: uncertainty trend}(b) demonstrates the changing trends of intersection-based indicator uncertainty and path reconstruction uncertainty regarding the number of UAVs for detection. The uncertainty of path reconstruction is calculated using Eq. \eqref{eq: path uncertainty}. It is obvious that both path reconstruction uncertainty and intersection-based indicator uncertainty decrease with the increase in UAVs for detection. 
The decreased amplitude of path reconstruction uncertainty is first larger than the intersection-based indicator uncertainty, while smaller than that with the increase of UAV-detected intersections, taking the case of 9-UAV deployment as a transition point.
It is sensible that when the number of deployed UAVs is less than 9, the marginal decrease of path reconstruction uncertainty is mainly due to the decreased number of unobserved paths, $\mathcal{K}-\mathcal{K}_o-\mathcal{K}_{cv}$, while the marginal decrease of path reconstruction uncertainty when the number of deployed UAVs is more than 9 UAVs is attributed to the transition of $u_{path}$ between $\mathcal{K}_{o,cv}$ and $\mathcal{K}_o-\mathcal{K}_{o,cv}$. By contrast, the relatively stable decreasing trend of intersection-based indicator uncertainty implies that the uncertainties of queue length and arrival rate are more focused on the intersection level, while the decrease of path reconstruction uncertainty is more affected by the correlated observability of multiple UAVs, considering the topology characteristics of paths within a network.

\begin{table*}[htbp]
\centering
\caption{\rmfamily\small Sensitivity analysis of weight coefficients.}
\label{tab: Sensitivity of weight coefficients}
\renewcommand\arraystretch{1.1}
\begin{tabular}{ccp{4.5cm}}
\hline
\textbf{\makecell[cc]{\rmfamily\small Weight coefficients}} & \textbf{\rmfamily\small Sum of sub-objective function value $F_{path} + F_{arrival} + F_{queue}$} \\
\hline
\rmfamily\small $w_1:w_2:w_3 = 6.5:1:1$ & \makecell[cc]{\rmfamily\small 2767} \\
\hline
\rmfamily\small $w_1:w_2:w_3 = 13:1:1$ & \makecell[cc]{\rmfamily\small 2767} \\
\hline
\rmfamily\small $w_1:w_2:w_3 = 26:1:1$ & \makecell[cc]{\rmfamily\small 2761} \\
\hline
\rmfamily\small $w_1:w_2:w_3 = 26:2:2$ & \makecell[cc]{\rmfamily\small 3011} \\
\hline
\rmfamily\small $w_1:w_2:w_3 = 26:4:4$ & \makecell[cc]{\rmfamily\small 3011} \\
\hline
\end{tabular}
\end{table*}

\subsubsection{Optimal UAV location and fleet size analysis}

Based on the results in Figure \ref{fig: uncertainty trend}, the optimal UAV location schemes under different UAV numbers are shown in Table \ref{tab: Optimal UAV deployment scheme}. 
When the coverage of UAVs is less than 50\%, the location of an added UAV is mostly the intersection with comparable traffic flow like Int. 1, 2, 3, 4 and 11, or the intersection with comparable passed path number like Int. 16 and 17, or the intersection featuring significant topological betweeness with neighboring intersection with comparable flow or passed path number. Taking the 9-UAV case as an example, which has the largest feasible space in this case study as the size of solution domain is the combinatorial number of selecting $\mathcal{I}_{uav}$ from $\mathcal{I}$, the optimal UAV deployment scheme is Int. 1, 2, 4, 5, 10, 11, 13, 16, 17, which are exactly the intersections with the significant flow and path number. It is noted that Int. 3 can be partially observed with its two links observable by UAVs, thus it is sensible that it is not selected. Considering the data provided by loop detectors at South Fuzhou Road, the intersections of South Fuzhou Road are not the priority for UAV deployment when the coverage of UAV is not high, which also implies the selection criteria of the proposed method to reduce the network-wide uncertainty.

\begin{table*}[htbp]
\centering
\caption{\rmfamily\small Optimal UAV deployment schemes given different detector numbers.}
\label{tab: Optimal UAV deployment scheme}
\renewcommand\arraystretch{1.1}
\begin{tabular}{ccp{4.5cm}}
\hline
\textbf{\makecell[cc]{\rmfamily\small UAV number}} & \textbf{\rmfamily\small UAV location scheme} \\
\hline
\rmfamily\small 1 & \makecell[cc]{\rmfamily\small 12}  \\
\hline
\rmfamily\small 2 & \makecell[cc]{\rmfamily\small 17, 18} \\
\hline
\rmfamily\small 3 & \makecell[cc]{\rmfamily\small 13, 16, 17} \\
\hline
\rmfamily\small 4 & \makecell[cc]{\rmfamily\small 4, 5, 16, 17} \\
\hline
\rmfamily\small 5 & \makecell[cc]{\rmfamily\small 4, 5, 13, 16, 17} \\
\hline
\rmfamily\small 6 & \makecell[cc]{\rmfamily\small 1, 4, 5, 13, 16, 17} \\
\hline
\rmfamily\small 7 & \makecell[cc]{\rmfamily\small 1, 2, 4, 10, 13, 16, 17} \\
\hline
\rmfamily\small 8 & \makecell[cc]{\rmfamily\small 1, 2, 4, 5, 11, 13, 16, 17} \\
\hline
\rmfamily\small 9 & \makecell[cc]{\rmfamily\small 1, 2, 4, 5, 10, 11, 13, 16, 17} \\
\hline
\rmfamily\small 10 & \makecell[cc]{\rmfamily\small 1, 2, 4, 5, 6, 10, 11, 13, 16, 17} \\
\hline
\rmfamily\small 11 & \makecell[cc]{\rmfamily\small 1, 2, 4, 5, 6, 9, 10, 11, 13, 16, 17} \\
\hline
\rmfamily\small 12 & \makecell[cc]{\rmfamily\small 1, 2, 4, 5, 6, 10, 11, 13, 14, 16, 17, 18} \\
\hline
\rmfamily\small 13 & \makecell[cc]{\rmfamily\small 1, 2, 3, 4, 5, 6, 10, 11, 13, 14, 16, 17, 18} \\
\hline
\rmfamily\small 14 & \makecell[cc]{\rmfamily\small 1, 2, 3, 4, 5, 6, 10, 11, 13, 14, 15, 16, 17, 18} \\
\hline
\rmfamily\small 15 & \makecell[cc]{\rmfamily\small 1, 2, 3, 4, 5, 6, 10, 11, 12, 13, 14, 15, 16, 17, 18} \\
\hline
\rmfamily\small 16 & \makecell[cl]{\rmfamily\small 1, 2, 3, 4, 5, 6, 7, 8, 10, 11, 13, 14, 15, 16, 17, 18} \\
\hline
\rmfamily\small 17 & \makecell[cl]{\rmfamily\small 1, 2, 3, 4, 5, 6, 7, 8, 10, 11, 12, 13, 14, 15, 16, 17, 18} \\
\hline
\end{tabular}
\end{table*}

\begin{figure}[!htbp]
\centering
\includegraphics[width=6 in]{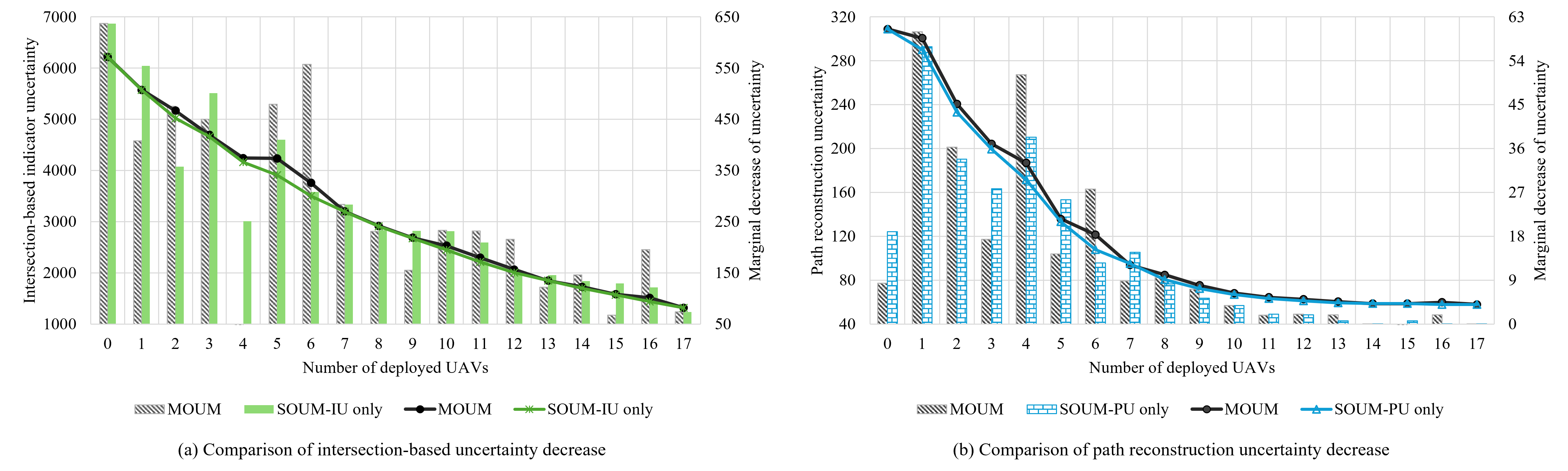}
\caption{\rmfamily\small Comparison of uncertainty decrease of different single-objective benchmarks }
\label{fig: uncertainty trend of SOUM}
\end{figure}

Further, the proposed uncertainty minimal method aimed at single objectives, intersection-based uncertainty (IU), and path reconstruction uncertainty (PU), is tested to demonstrate the sensitivity difference of uncertainty decrease under different optimization objectives. Correspondingly, the decrease of IU and PU in Figure \ref{fig: uncertainty trend} is used for comparison in Figure \ref{fig: uncertainty trend of SOUM}. Though the proposed MOUM demonstrates a slightly larger uncertainty as compared with both benchmarks of SOUM in some cases, the changing trends of uncertainty decrease of the two SOUM benchmarks and the proposed MOUM are similar. Taking the 5-UAV case in Figure \ref{fig: uncertainty trend of SOUM} (a) for example, the optimal UAV deployment scheme obtained by SOUM-IU is Int. 1, 5, 13, 17, 18, while the optimal scheme by the proposed MOUM is Int. 4, 5, 13, 16, 17, as shown in Table \ref{tab: Optimal UAV deployment scheme}. The different UAV locations at Int. 1 and 18 have a larger intersection-based uncertainty than Int. 4 and 16, as shown in Figure \ref{fig: uncertainty trend} (a), which accounts for the slightly larger uncertainty of the scheme of MOUM. Such a gap can be explained through the larger number of passed paths of Int. 4 and 16 in Fig. \ref{fig: shinan_network}(b), because the path reconstruction uncertainty also affects the selection of UAVs as one of the objectives. Such a difference to some degree validates the effectiveness of adopting a multi-objective model and the ability of the solution algorithm in optimum searching.

\begin{figure}[!htbp]
\centering
\includegraphics[width = 3.4 in]{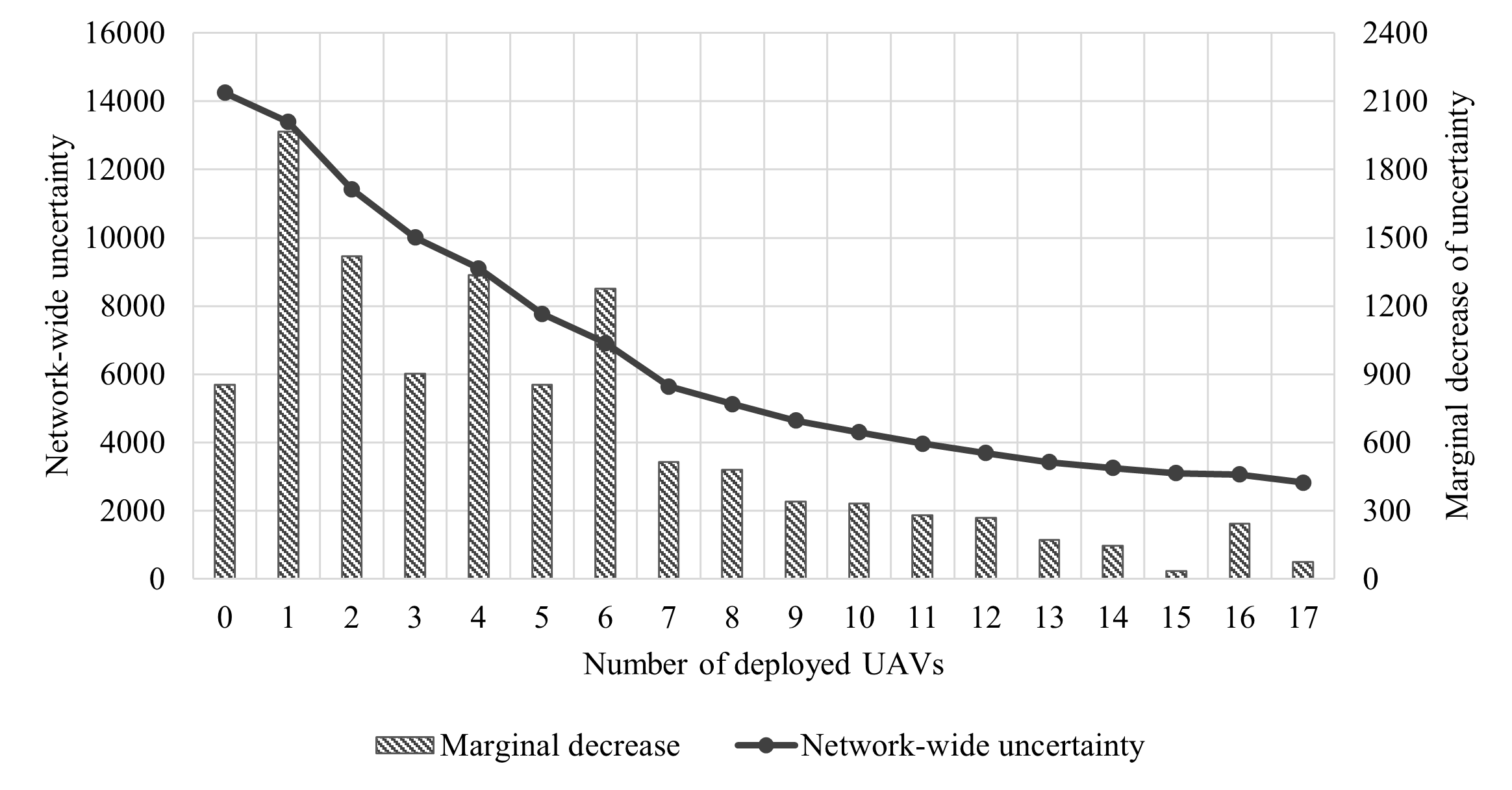}
\caption{\rmfamily\small Sensitivity of network-wide uncertainty to the number of deployed UAVs  }
\label{fig: optimal UAV fleet size}
\end{figure}

It is noted that the bar plots of marginal uncertainty for the case of 0 UAV refer to the uncertainty decrease from 0 UAV to 1 UAV, and so on. Considering the marginal uncertainty decrease, the case of a 6-UAV fleet size can be regarded as a critical point, as the uncertainty decrease becomes flatter when the number of UAVs is larger than 7. Thus, the optimal UAV fleet size is 7 in this case study, which is also verified by the changing trend of network-wide uncertainty $Z$ as shown in Figure \ref{fig: optimal UAV fleet size}. Further comparison on the effectiveness of the optimal UAV location schemes is introduced in Subsection \ref{sec:horizontal comparison}.

\subsection{Horizontal comparison} \label{sec:horizontal comparison}
As a variant of the set coverage problem, the proposed UAV location optimization problem mostly adopts the flow coverage maximum as the objective in existing studies. Considering the ultimate goal of network-wide traffic state monitoring, a horizontal comparison is conducted with a flow coverage maximum (FCM) model and the single-objective benchmark regarding the downstream application of queue length estimation, arrival rate estimation, and path flow estimation. 

The selected FCM model \citep{GENTILI2012227} is formulated as Eq. \eqref{eq: coverage maximum} - \eqref{eq: binary variable}, where $x_a$ is a binary variable for each link $a \in \mathcal{A}$ which is equal to 1 if link $a$ is observed and otherwise. $v_a$ denotes the vehicle flow volume of link $a$. Solved with greedy algorithm and classic genetic algorithm (GA), the FCM benchmark was evaluated with $N_{uav}$ ranging from 1 to 17, with three parallel tests for each certain UAV number. 

\begin{align}
    max \quad \sum_{a \in \mathcal{A}}  v_{a} x_{a} \label{eq: coverage maximum} \\
    s.t. \quad \sum_{a \in \mathcal{A}} x_{a} \leq N_{uav} \label{eq: limited UAV} \\
    x_a \in \{0,1\} \label{eq: binary variable}
\end{align}

\begin{figure}[!htbp]
\centering
\includegraphics[width=5.8 in]{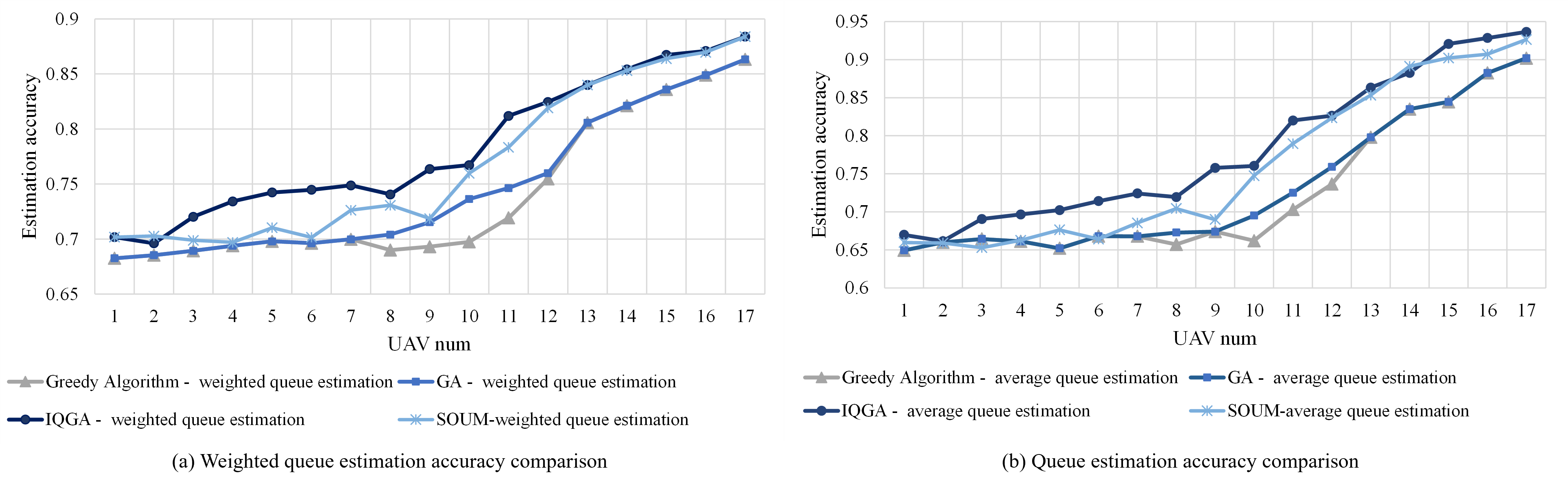}
\caption{\rmfamily\small Queue estimation comparison between different methods.}
\label{fig: queue estimation comparison}
\end{figure}

\begin{figure}[!htbp]
\centering
\includegraphics[width=5.8 in]{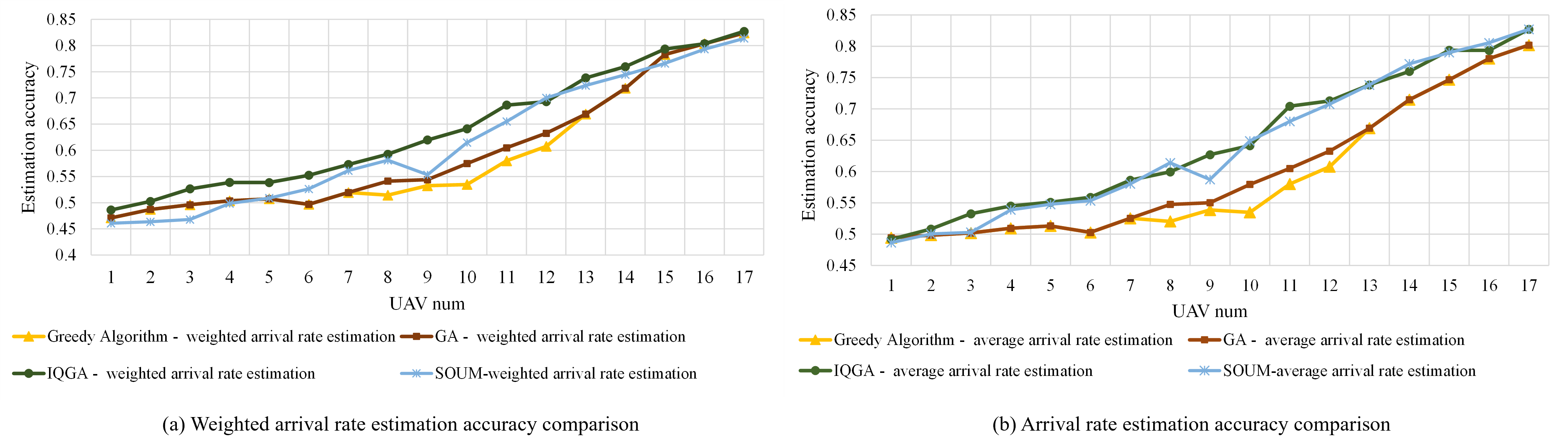}
\caption{\rmfamily\small Arrival rate estimation comparison between different methods.}
\label{fig: arrival rate estimation comparison}
\end{figure}

\begin{figure}[!htbp]
\centering
\includegraphics[width=5.8 in]{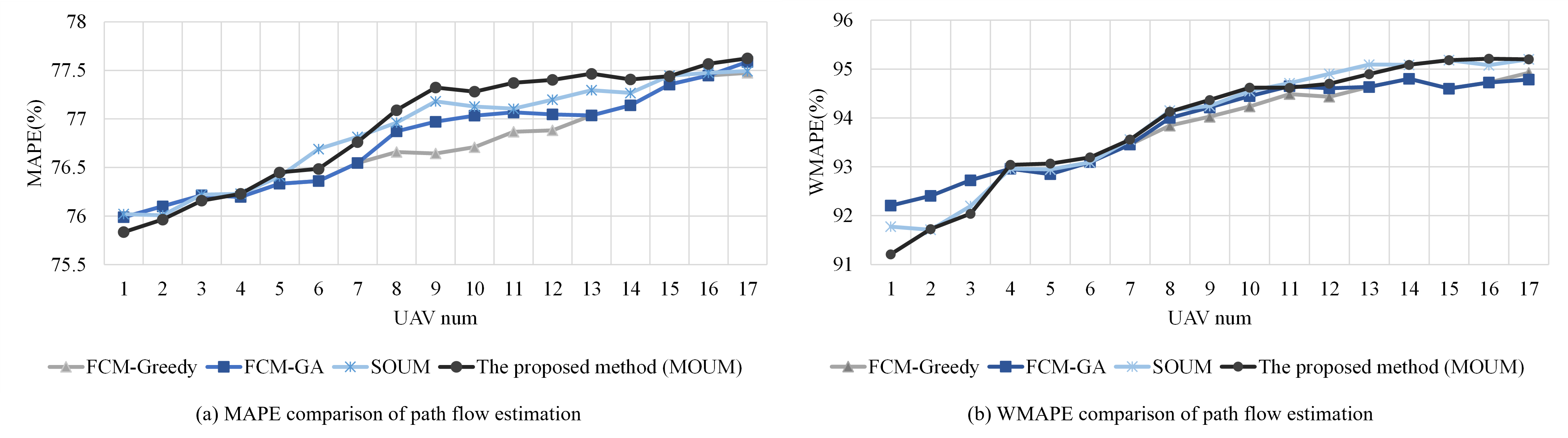}
\caption{\rmfamily\small Path flow estimation comparison between different methods.}
\label{fig: path flow estimation comparison}
\end{figure}

As for the selected methods for downstream traffic state estimation to evaluate the effectiveness of UAV location optimization, one cycle-based queue length estimation method \citep{2017Real}, one cycle-based arrival rate estimation method \citep{tan2025} and one path flow estimation method \citep{yao2022path} were used as the benchmarks.
Based on the output optimal UAV location schemes of the proposed MOUM method, two SOUM baselines and the FCM benchmark, the multi-source data were input for traffic state estimation. Using the ground truth extracted from the VISSIM simulation run, the comparison of downstream traffic state estimation was shown in Fig. \ref{fig: queue estimation comparison}, Fig. \ref{fig: arrival rate estimation comparison}, and Fig. \ref{fig: path flow estimation comparison}. It is noted that the SOUM-IU baseline is used for comparison of intersection-based indicator estimation in Fig. \ref{fig: queue estimation comparison} and Fig. \ref{fig: arrival rate estimation comparison}, including queue length and arrival rate, while the SOUM-PU baseline is used for path flow estimation in Fig. \ref{fig: path flow estimation comparison}. To quantify the estimation accuracy, the following evaluation indicators are chosen to plot such figures. 
\begin{itemize}
    \item \textbf{Weighted queue estimation accuracy}. The flow-weighted average of the Mean Absolute Percentage Error (MAPE) of arrival rate estimation, which is calculated by \\$1-\frac{1}{\sum_{m \in \mathcal{M}_i} Q_{m}}{\sum_{m \in \mathcal{M}_i} (\frac{Q_{m}} {card(\mathcal{J}_{m})} * \sum_{j \in \mathcal{J}_{m}} \sum_{t \in \{1,2, \dots, C_j/\Delta\}} \frac{|\hat{A}_{t}-A_{t}|}{A_{t}})}$.
    \item \textbf{Queue estimation accuracy}. The average of the Mean Absolute Percentage Error (MAPE) of BoQ estimation, which is calculated by $1-\frac{1}{ card(\mathcal{M}_i)}{\sum_{m \in \mathcal{M}_i} (\frac{1}{card(\mathcal{J}_{m})} * \sum_{j \in \mathcal{J}_{m}} \frac{|\hat{L}_{j}-L_{j}|}{L_{j}})}$.
    \item \textbf{Weighted arrival rate estimation accuracy}. The flow-weighted average of the Mean Absolute Percentage Error (MAPE) of arrival rate estimation, which is calculated by \\$1-\frac{1}{\sum_{m \in \mathcal{M}_i} Q_{m}}{\sum_{m \in \mathcal{M}_i} (\frac{Q_{m}} {card(\mathcal{J}_{m})} * \sum_{j \in \mathcal{J}_{m}} \sum_{t \in \{1,2, \dots, C_j/\Delta\}} \frac{|\hat{A}_{t}-A_{t}|}{A_{t}})}$.
    \item \textbf{Arrival rate estimation accuracy}. The average of the Mean Absolute Percentage Error (MAPE) of BoQ estimation, which is calculated by $1-\frac{1}{card(\mathcal{M}_i)} {\sum_{m \in \mathcal{M}_i} (\frac{1}{card(\mathcal{J}_{m})} * \sum_{j \in \mathcal{J}_{m}} \sum_{t \in \{1,2, \dots, C_j/\Delta\}} \frac{|\hat{A}_{t}-A_{t}|}{A_{t}})}$.
    \item \textbf{Weighted path flow estimation accuracy}. The Weighted Mean Absolute Percentage Error of path flow estimation, which is calculated by $1-\frac{\sum_{k \in \mathcal{K}}|\hat{Q}_{k}-Q_{k}|}{\sum_{k \in \mathcal{K}} Q_k}$.
    \item \textbf{Path flow estimation accuracy}. The average of the Mean Absolute Percentage Error (MAPE) of path flow estimation, which is calculated by $1-\frac{1}{card(\mathcal{K})} \sum_{k \in \mathcal{K}} \frac{|\hat{Q}_{k}-Q_{k}|}{Q_k}$.
    \item \textbf{Path coverage}. The number of paths observed by the multi-source detector, including CV, loop and UAV.
    \item \textbf{Flow coverage}. The sum of movement flows observed by the multi-source detector, including CV, loop, and UAV.
\end{itemize}

It is noted that $Q_{m}$ denotes the flow of movement $m$. $card()$ is the operator of obtaining the size of a set. $\hat{L_{j}}, L_{j}$ denote the estimation and the ground truth of cycle $j$, respectively. $\hat{A}_{t}, A_{t}$ denote the estimation and the ground truth of unit time step $t$, respectively. $C_j$ denotes the cycle length of cycle $j$. $\Delta$ is the duration of a unit time step. $\hat{Q}_{k}, Q_{k}$ denote the estimation and the ground truth of path flow of path $k$,  respectively. 

\subsubsection{Intersection-based indicator estimation comparison}

All methods present an increasing trend of the estimation accuracy with the increase of deployed UAV number, and the proposed method outperforms the selected baselines in nearly all the cases covering different UAV fleet size cases.
For the FCM baseline, the estimation accuracies of the solution obtained by GA are better than those of the solution obtained by gthe reedy algorithm when the number of deployed UAVs ranges from 8 to 12 by about 0.89\% $\sim$ 1.42\%, which shows the ability of the solution searching of the metaheuristic algorithm in medium and large-scale problems. As compared with the solution of GA, the estimation accuracies of the proposed method are better by about 7.23\% in arrival rate estimation and 5.02\% in queue length estimation, and gaps are higher in flow-weighted arrival rate and queue length estimation by about 8.04\% and 8.08\%, respectively. The increased advantages of the weighted accuracy imply that UAV deployment aimed at the objective of maximum flow coverage is not necessarily a guarantee of the estimation accuracy of the network-wide traffic monitoring. 

As for the SOUM baseline, the estimation accuracies of queue length and arrival rate are similar to those of the proposed MOUM method when the number of deployed UAVs is larger than 12. Except for the average arrival rate estimation, the gap between the proposed method and the SOUM baseline is obvious when the number of deployed UAVs is smaller than 12, which can reach 6.79 \%. Especially for the 9-UAV case which has the largest solution space, the superiority of the proposed method shows the effectiveness of traffic state estimation-oriented UAV location optimization.

\subsubsection{Path flow estimation comparison}

Similar to the estimation of intersection-based indicators, the path flow estimation accuracy increases with the growth of the UAV fleet size. Except for the cases of less than 4 UAVs, the estimation accuracy is larger than the FCM baseline and the SOUM baseline in most cases, which may be due to the preference for increased observability by paired UAVs in the model formulation as well as the solution algorithm. Though the gaps of path flow estimation accuracy are slightly 0.132\% and 0.265\% for WMAPE and MAPE, respectively, which is far smaller than its superiority of intersection-based indicator estimation, the proposed method does show its advantages in downstream traffic state estimation.

\begin{figure}[!htbp]
\centering
\includegraphics[width=3.6 in]{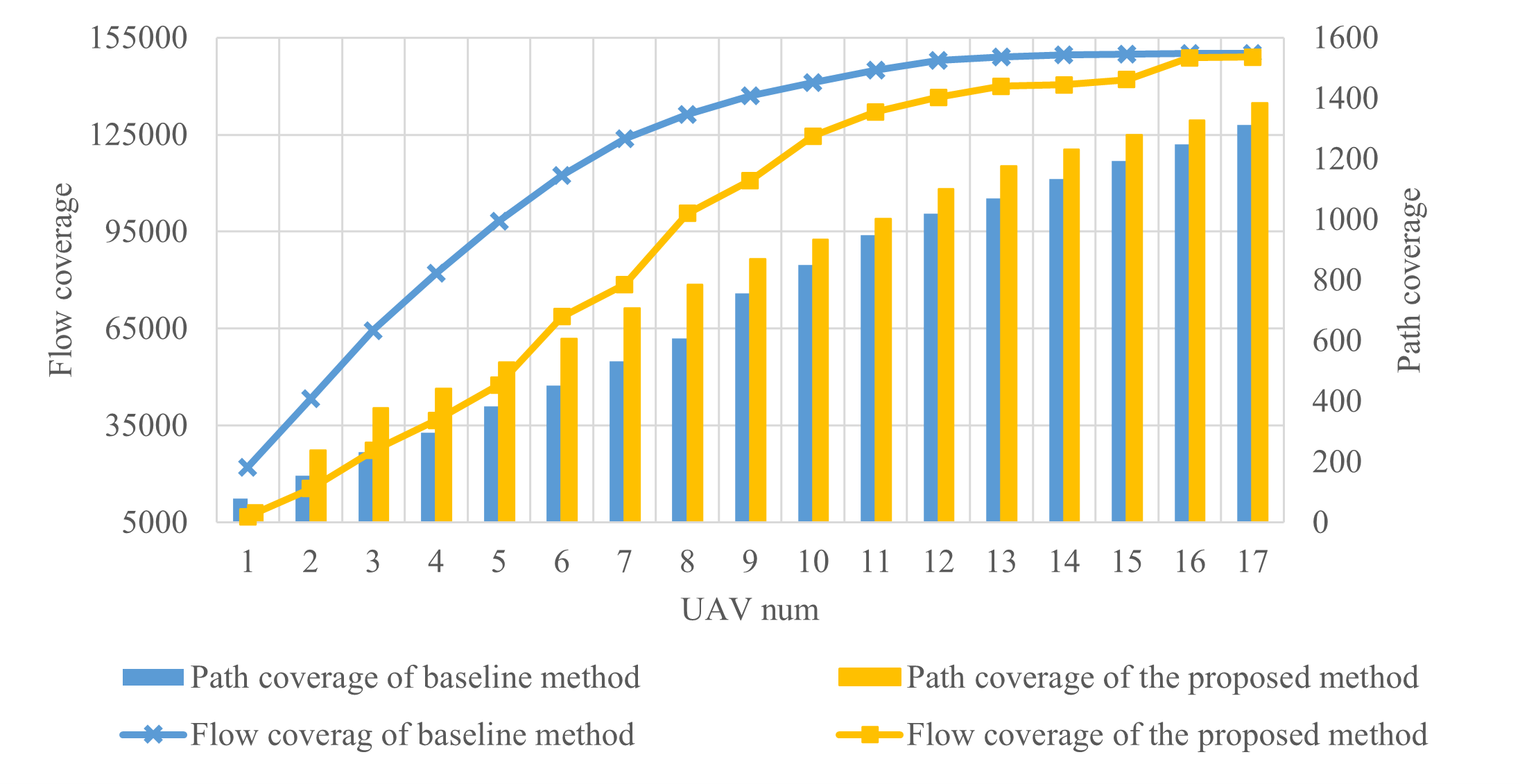}
\caption{\rmfamily\small Comparison of flow and path coverage between different methods.}
\label{fig: flow coverage comparison}
\end{figure}

The reason to the effectiveness of the proposed method in the downstream traffic state estimation application can also be implied through the comparison of flow and path coverage as shown in Figure \ref{fig: flow coverage comparison}, which plots the comparison between the solution obtained by the FCM-GA baseline and the proposed method. Although the FCM baseline shows superiority of link flow coverage over the proposed method by about 27.29\%, the coverage of passed paths of the proposed method is larger than that of the FCM baseline by about 20.47\%. Considering the downstream traffic state estimation application, a high path coverage actually reveals more information about the path flows, which constitute the movement flows for intersection-based indicator estimation. In this way, it can be concluded that a high flow coverage does not necessarily relate to a high estimation accuracy for traffic state monitoring; instead, it is the critical information of path choice and arrival pattern showing the joint effect of both traffic demand and supply that matters. Besides, the gap between both the flow and path coverage presents a decreasing trend with the increasing number of UAVs to be deployed, which, to some degree, justifies the significance of a reasonable location scheme to fully exert the marginal benefits when the number of available UAVs is limited for network-wide traffic state monitoring.

\subsection{Solution algorithm performance}
\begin{figure}[!htbp]
\centering
\includegraphics[width=6 in]{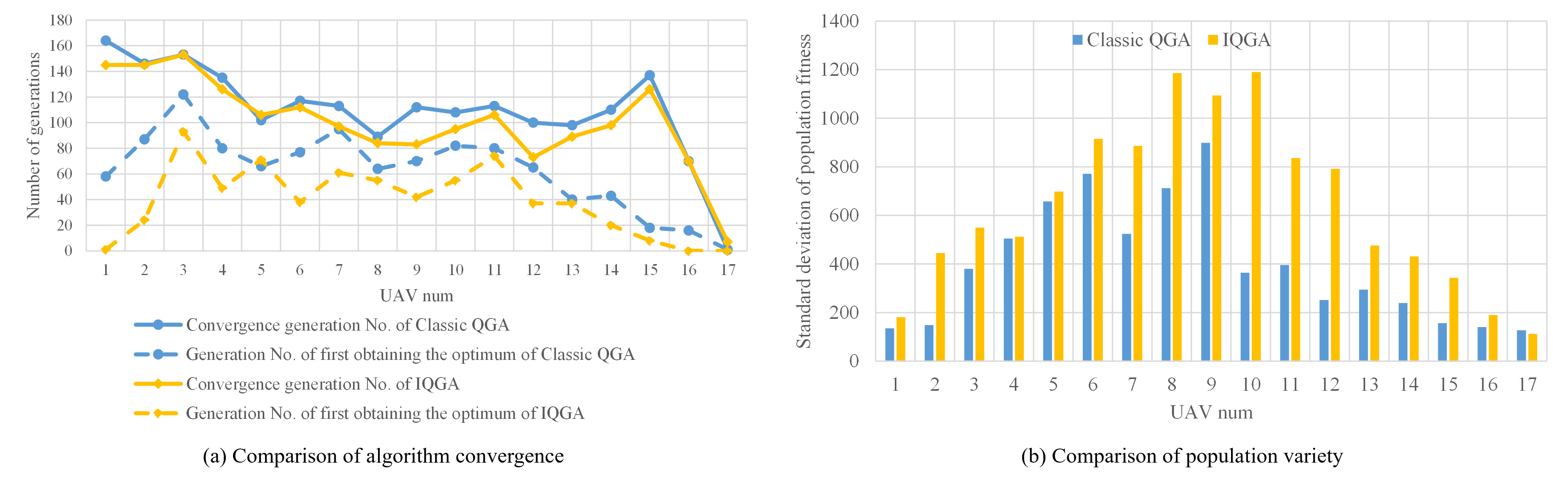}
\caption{\rmfamily\small Algorithm improvement of IQGA on algorithm convergence and solution variety.}
\label{fig: IQGA improvement}
\end{figure}

In order to evaluate the effectiveness of the proposed IQGA, a horizontal comparison is conducted with classic QGA without the two added operators, while the parameters are the same as IQGA. As IQGA is improved with the deduplication operator and fine-tuning operator for better exploration and exploitation, the generation number to obtain the optimum and reach convergence, as well as the variance of population fitness, are compared as shown in Figure \ref{fig: IQGA improvement}. 
The selected evaluation indicators are defined as below. 
\begin{itemize}
    \item \textbf{Convergence generation number}. The index of population evolution generation where the fitness function starts to converge.
    \item \textbf{Generation number of first obtaining the optimum}. The index of population evolution generation where the fitness function first obtains the optimum.
    \item \textbf{Standard deviation of population fitness}. The standard deviation of fitness of the whole population for each evolution generation.
\end{itemize}

It can be seen that IQGA can find the optimum with fewer generations, meanwhile realizing faster convergence, under the joint effect of two added operators. Besides, the standard deviation of the population fitness of IQGA is higher than that of classic QGA by an average of 77.69\%, especially when the number of deployed UAVs ranges from 7 to 12, which have relatively larger solution spaces. 

\begin{figure}[!htbp]
\centering
\includegraphics[width=6 in]{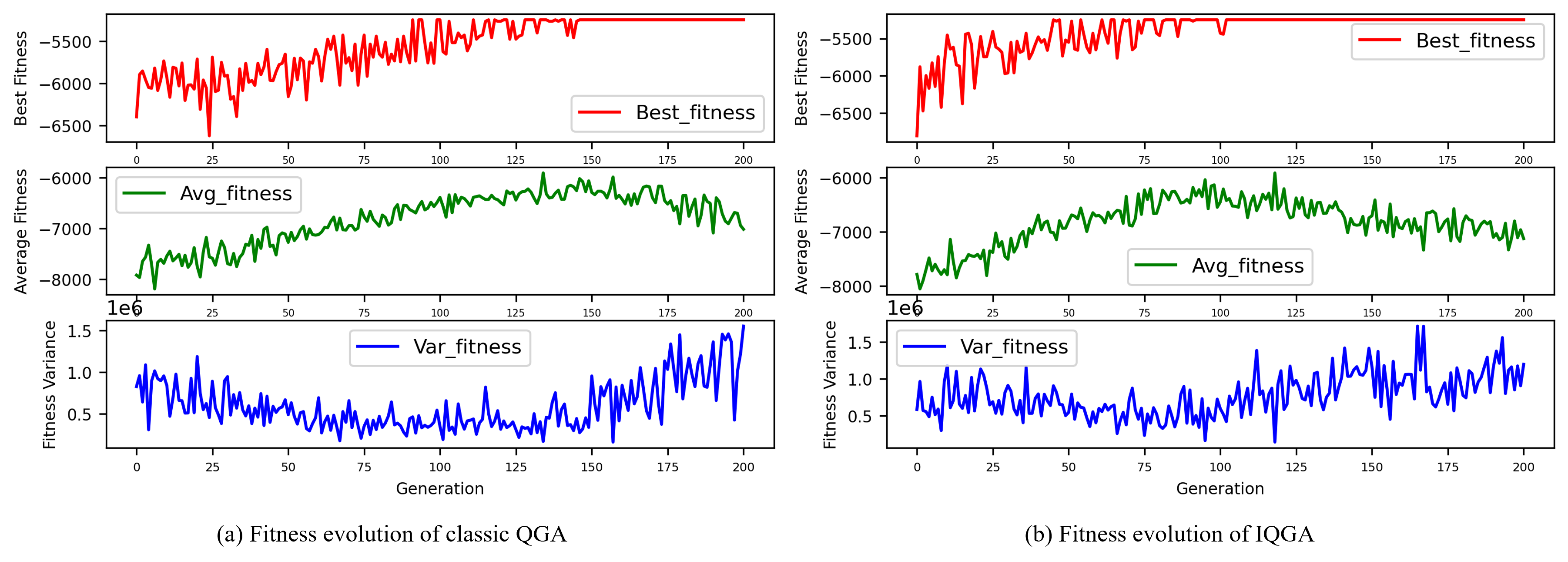}
\caption{\rmfamily\small Algorithm performance comparison under the case of 9 deployed UAVs.}
\label{fig: 9-UAV case comparison}
\end{figure}

Taking the 9-UAV case as an example, the mean value, optimal value, and variance of the fitness function of the population for each evolution generation are plotted in Figure \ref{fig: 9-UAV case comparison}. It can be seen that the population evolution process of IQGA converges earlier to the best fitness, while the variance of fitness is relatively larger throughout all the generations.

\begin{figure}[!htbp]
\centering
\includegraphics[width=3.2 in]{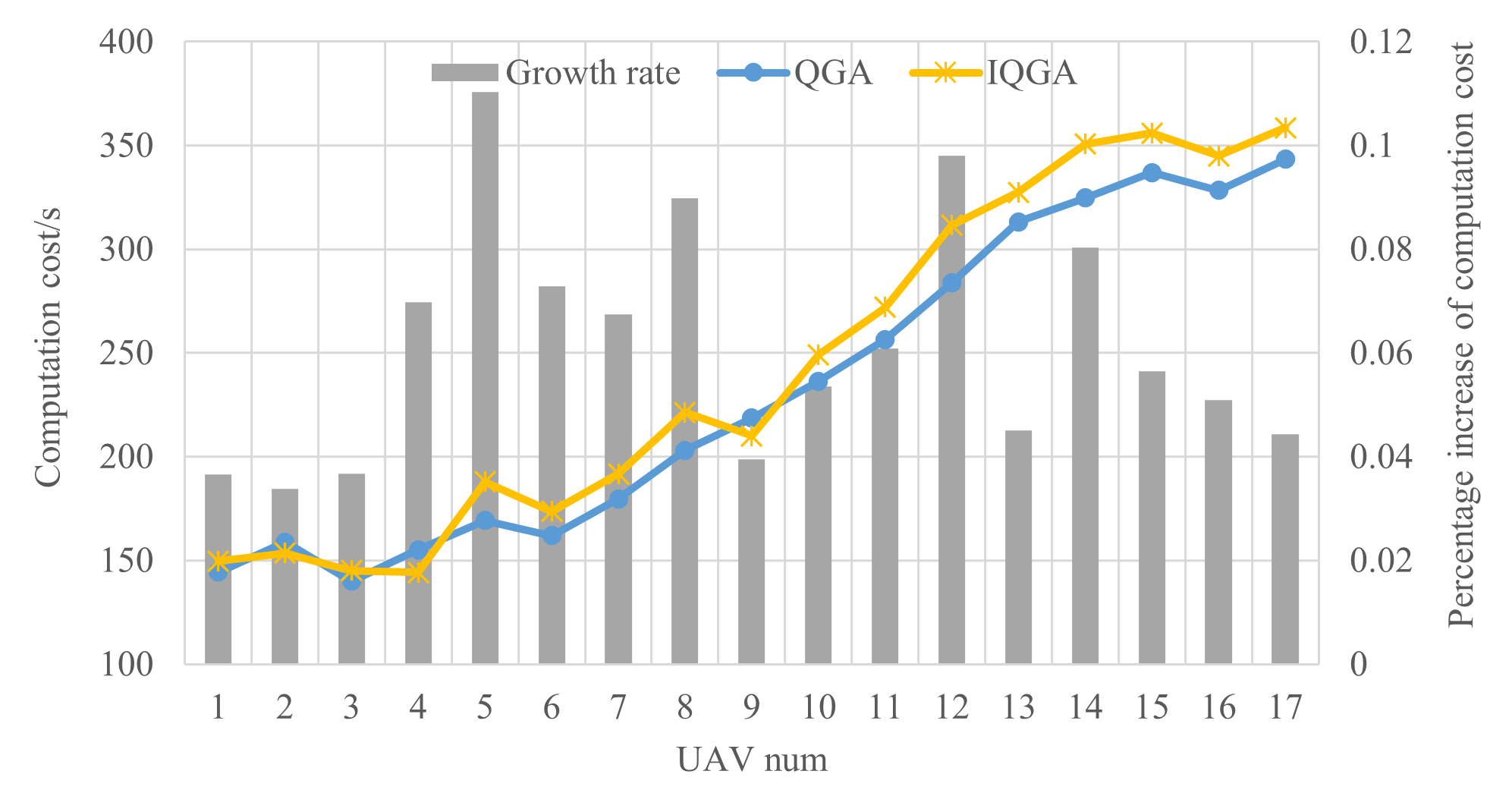}
\caption{\rmfamily\small Computation cost comparison.}
\label{fig: computation cost}
\end{figure}

Further, the computation cost of QGA and IQGA under the 9-UAV case is shown in Figure \ref{fig: computation cost}. It is noted that the percentage difference of computation cost of the secondary vertical axis refers to the percentage difference of computation time between different solution algorithms, which is calculated by $\frac{|{T}_{b}-T'|}{T_b}$, where ${T}_{b}$ represents the computation time of the algorithm used as the baseline, while $T'$ is the computation time of the algorithm used for comparison.
Though the computation cost increase caused by the added operators ranges from 3.38\% to 11.03\%, the computation cost of IQGA can actually be reduced by setting a smaller generation number, considering the faster convergence speed with guarantee of solution quality.

\section{Conclusion and Future Work}

This study explores the optimal UAV deployment problem in a signalized road network in conjunction with ground sensor data to achieve a more reliable estimation of network-wide traffic states, including path flow, arrival rates and queue lengths. Specifically, the optimal UAV deployment problem is transformed into an optimization problem to minimize the uncertainty of network-wide traffic states, where an indicator of feasible domain size and the corresponding calculation approaches are proposed to measure the uncertainty of different network-wide traffic states. 

Evaluation results of a simulation case study have demonstrated that the optimal UAV location scheme obtained by the proposed method can reach an improvement of 7.23\%, 5.02\% and 0.265\% in terms of the estimation accuracy of arrival rate, queue length, and path flow, respectively. The proposed IQGA is also shown to be faster in solution convergence than the classic QGA by about 9.22\%. The effectiveness of modeling different aspects of traffic state uncertainty as the objective of UAV location and the ability of IQGA in searching the optimum is thus well proven through its superiority over baselines in the downstream traffic state estimation application. 

Future research includes i) sensitivity analysis on factors like the CV penetration rate as well as parameters of QGA; ii) further testing of the enhancement of UAV deployment on different traffic state estimation methods, like travel time, OD flow, etc.; iii) extension of UAV location to UAV routing aimed at dynamic monitoring tasks considering all-day traffic evolution.

\bibliographystyle{cas-model2-names}

\bibliography{Arxiv_version}

\end{document}